\documentstyle{amsppt}
\magnification 1050
\vsize 22.5truecm
\NoBlackBoxes
\NoRunningHeads
\topmatter
\title Hopf Hypersurfaces of Low Type in Non-flat Complex Space Forms
\endtitle
\author Ivko Dimitri\'c
\endauthor
\abstract We classify Hopf hypersurfaces of non-flat complex space forms $\Bbb CP^m(4)$ and $\Bbb CH^m(-4),\,$ denoted jointly by $\Bbb CQ^m(4c)$,
that are of 2-type in the
sense of B. Y. Chen, via the embedding into a suitable (pseudo) Euclidean space of Hermitian matrices by projection operators.
This complements and extends earlier classifications by Martinez-Ros (minimal case)  and Udagawa (CMC case), who studied only hypersurfaces of
$\Bbb CP^m$ and assumed them to have constant mean curvature instead of being Hopf. Moreover, we rectify some claims in Udagawa's paper to
give a complete classification of constant-mean-curvature-hypersurfaces of 2-type. We also derive a certain characterization of CMC
Hopf hypersurfaces which are of 3-type and mass-symmetric in a naturally-defined  hyperquadric containing the image of $\Bbb CQ^m(4c)$
via these embeddings. The classification of such hypersurfaces is done in $\Bbb CQ^2(4c), $ under an additional assumption in the hyperbolic
case that the mean curvature is not equal to $ \pm 2/3.$ In the process we show that every standard example of class $B$ in $\Bbb CQ^m(4c)$
is mass-symmetric and we determine its Chen-type.
\endabstract
\footnote[]{The work partially supported by the Eberly Family Science Endowment Fund, PSU Fayette.}
\address Mathematics Department, Penn State
University Fayette, The Eberly Campus, P.O.  Box 519, Uniontown, PA 15401, USA. \endaddress
\email ivko\@psu.edu \endemail
 \endtopmatter \document

\define\k{\varkappa}
\define\s{\sigma}
\define\n{\nabla f}

\define\cxf{\Bbb CQ^m}
\define\lp{\Delta}
\define\la{\langle}
\define\ra{\rangle}
\define\tx{\tilde x}

\head 1. Introduction
\endhead

\par The study of finite-type submanifolds of Euclidean and pseudo-Euclidean spaces
has been an area of flourishing research initiated by
B. Y. Chen in the 1980s [9]. Many geometers contributed to the theory and quite a number of important and interesting results coming from that 
study have been
obtained on sharp eigenvalue estimates and characterizations of certain submanifolds by eigenvalue equalities [10].
A Riemannian $n-$manifold $M^n$ isometrically immersed into a Euclidean or pseudo-Euclidean space by $x: M^n \to E^N_{(K)} $  is said to be
of $k-$type (more precisely of Chen $k-$type) in $E^N_{(K)}$ if the position vector $x$ can be decomposed, up to a
translation by a constant vector $x_0$, into a sum of $k$ nonconstant
$E^N_{(K)}-$valued eigenfunctions of the Laplacian $\lp_M$ from different eigenspaces, viz.
$$ x = x_0 + x_{t_1} + \dots + x_{t_k}; \quad   x_0 = \text{const}, \quad
\lp x_{t_i} = \lambda_{t_i} x_{t_i}, \;\; i= 1,...,k, \tag 1 $$
where $x_{t_i} \neq \text{const}, \;\lambda_{t_i} \in \Bbb R$ are all different, and the Laplacian acts on a vector-valued function componentwise.
For a compact submanifold, the constant part $x_0$ is the center of mass and if $x$ immerses $M^n$ into a central hyperquadric
of a Euclidean or pseudo-Euclidean space the immersion is said to be mass-symmetric in that hyperquadric if $x_0$ coincides with the center
of the said hyperquadric. Moreover, decomposition (1) also makes sense for noncompact submanifolds, but $x_0$ may not be uniquely determined,
namely when one of the eigenvalues $\lambda_{t_i}$ above is zero. Such submanifolds are said to be of null $k$-type, and are, therefore, per 
definition mass-symmetric.
\par
The study of finite-type submanifolds therefore treats an interesting question: To what extent is the geometric structure of a
submanifold determined by a simple analytic information, that is, by the spectral resolution (1) of the immersion into
finitely many terms?  By placing a complex projective or a complex hyperbolic space into a suitable (pseudo) Euclidean
space of Hermitian matrices using the embedding $\Phi$ by projectors in the standard way (cf.  [31], [28], [29], [18], [15]), it is possible to study
submanifolds, in particular hypersurfaces, of a complex space form in
terms of finite-type property, where the immersion considered is the composite immersion with $\Phi$.  It is well-known that a 1-type submanifold is minimal in an appropriate hyperquadric of the ambient (pseudo) Euclidean space. 1-Type real hypersurfaces of a complex space
form  $\cxf$ were previously studied in [23], [18], and the present author subsequently  classified 1-type
submanifolds of these spaces of any dimension (see [14], [15]).
In particular, 1-type hypersurface in $\Bbb CP^m(4)$ is a geodesic hypersphere of radius
$r = \arctan\sqrt{2m + 1},$ which has an interesting stability property [23], [15].
Type-2 (also called bi-order) hypersurfaces in the complex projective space were studied
by Martinez and Ros [23] (the minimal case) and  Udagawa [36], who classified them under
the assumption that they have constant mean curvature. However, Udagawa's classification is incomplete and has some deficiencies which we 
rectify here. 
First, it was claimed without proof in [36, p.194] that there are no 2-type  hypersurfaces in $\Bbb CP^m$ among homogeneous examples of 
class $B.$ We find a counterexample to this claim, producing two such hypersurfaces. Second,  it was claimed in the same paper 
(pp. 192-193) that there are no geodesic hyperspheres (i.e. class-$A_1$ hypersurfaces) in $\Bbb CP^m$ which are mass-symmetric and of 2-type, 
whereas we prove that a geodesic hypersphere of radius $\cot^{-1}(1/\sqrt m)$ exactly has these properties. Because of these erroneous claims, 
all three theorems of [36] are deficient in one way or another. 
\par
K\"ahler submanifolds of $\Bbb CP^m(4)$ of
2-type were successfully studied and classified in works of Ros [29] and Udagawa [35], whereas Shen [30] produced
a classification of minimal surfaces (real dimension 2) in $\Bbb CP^m (4)$ of 2-type. On the other hand, there are only scant
results so far on 3-type submanifolds of complex space forms (see [33], [34]) and their further study is warranted.
An overview of the results on low-type submanifolds of projective and hyperbolic spaces via the immersion by projectors is presented in [16].
\par
 In this paper we further advance the
study of hypersurfaces of non-flat complex space forms (that is, of both complex
projective and complex hyperbolic space) which are of 2- or
3-type and produce some new classification results, with the starting (weaker)
assumption that the hypersurfaces possess some simple compatibility property between the
complex structure of the ambient space and the second fundamental form.  One of the most
studied kinds of hypersurfaces in complex space forms are the so-called Hopf hypersurfaces [3], [11], [25],
defined by the property that the (almost contact) structure vector $U:= -J\xi$,
where $\xi$ is the unit normal, is a principal curvature vector (i.e. proper for the shape operator).
Equivalently, they are defined by integral curves of
the structure vector field $U$ being geodesics and in $\Bbb CP^m$ they are realized as tubes about complex submanifolds when the corresponding
focal set has constant rank [8]. The above-mentioned examples of
2-type hypersurfaces studied in [23] and [36] are in fact certain homogeneous Hopf
hypersurfaces. One of our results is that a 2-type
Hopf hypersurface indeed has constant mean curvature, the key result towards their classification
given in Theorems 1 and 2.
K\"ahler submanifolds of 3-type in complex projective spaces are
studied in [33], [34], where some examples are given, including compact irreducible Hermitian symmetric submanifolds of degree 3.
In this paper we also undertake a study of 3-type Hopf hypersurfaces with constant mean curvature in non-Euclidean complex space forms, 
fulfilling the promise made in [13], based on the study of spherical hypersurfaces of constant mean curvature which are of 3-type via 
the second standard immersion of the unit sphere (see also [17]).
Along the way we obtain a generalization of
Nomizu-Smyth's formula for the trace Laplacian of the shape operator [26],
and Simons'-type formula for the Laplacian of the squared norm of the second
fundamental form, which may be useful in other contexts.
For the background and additional clarification of the  notation used in this article a reader
should consult [15]. Excellent references on the geometry of hypersurfaces of complex space
forms are [3], [4], [25], and a brief overview [5].
\vskip 0.3truecm

\head 2. The basic background and relevant formulas
\endhead

Let $\Bbb CQ^m(4c)$ denote $m-$dimensional non-flat model complex space form, that is either the complex projective space $\Bbb CP^m(4)$
or the complex hyperbolic space $\Bbb CH^m(-4)$ of constant holomorphic sectional curvature $4c \, (c = \pm 1).$ By using a particular (pseudo) Riemannian submersion one can construct $\Bbb CQ^m$ and its embedding $\Phi$ into a certain (pseudo) Euclidean space of matrices. Consider first Hermitian form $\Psi_c$ on $\Bbb C^{m+1}$ given by $\Psi_c (z, w) = c \bar z_0 w_0 + \sum_{j=1}^m \bar z_j w_j, \; z, w \in \Bbb C^{m+1}$
with the associated (pseudo) Riemannian metric $g_c = \text{Re}\,\Psi_c$ and the quadric hypersurface
$N^{2m+1}:= \{ z \in \Bbb C^{m+1}\, | \, \Psi_c (z, z) = c\}.$ When $c = 1, \; N^{2m+1}$ is the ordinary hypersphere $S^{2m+1}$ of 
$\Bbb C^{m+1} = \Bbb R^{2m+2}$ and when $c = -1, \; N^{2m+1}$ is the anti - de Sitter space $H_1^{2m+1}$ in $\Bbb C^{m+1}_1.$ The orbit 
space under the natural action of the circle group $S^1$ on $N^{2m+1}$ defines $\Bbb CQ^m(4c),$ which is then the base space of a 
(pseudo) Riemannian submersion with totally geodesic fibers. The standard embedding $\Phi$ into the set of $\Psi-$Hermitian matrices 
$H^{(1)}(m+1)$ is achieved by identifying  a point, that is a complex line (or a time-like complex line in the hyperbolic case) with the 
projection operator onto it. Then one gets the following matrix representation of $\Phi$ at a point $p = [z],$ where $z = (z_j) \in 
N^{2m+1}\subset \Bbb C^{m+1}_{(1)}$
$$ \Phi ([z]) = \pmatrix |z_0|^2 & c z_0 \bar z_1 & \cdots & c z_0\bar z_m \\ z_1 \bar z_0 & c |z_1|^2  & \cdots & c z_1\bar z_m\\
\vdots & \vdots & \ddots & \vdots \\ z_m \bar z_0 & c z_m \bar z_1 & \cdots & c |z_m|^2
\endpmatrix . \tag 2
$$
The second fundamental form $\s $ of this embedding is parallel and the image $\Phi (\Bbb CQ^m)$ of the space form is contained in the hyperquadric of $H^{(1)}(m+1)$ centered at $I/(m+1)$ and defined by the equation
$$\langle P - I/(m+1),\, P - I/(m+1)\rangle = \frac{cm}{2(m+1)},
$$
where $I$ denotes the $(m+1)\times (m+1)$ identity matrix.
For the fundamental properties of the embedding $\Phi$ see [31], [18], [28], [15].
\par
If now $x: M^n \to \Bbb CQ^m (4c)$ is an isometric immersion of a Riemannian $n-$manifold as a real hypersurface of a complex space form ($n = 2m-1$) then we have the associated composite immersion $\tx = \Phi \circ x,$ which realizes $M$ as a submanifold of the (pseudo) Euclidean space
$E^N_{(K)}:= H^{(1)}(m+1),$ equipped with
the usual trace metric $\la A, B \ra = \frac c 2 \text{tr}\, (AB) $. In this notation the subscripts and superscripts in parenthesis are present
only in relation to $\Bbb CH^m,$ so that the superscript 1 in $H^{(1)}(m+1)$ is optional and appears only in the hyperbolic case, since the
construction of the embedding is based on the form $\Psi$ in $\Bbb C_1^{m+1}$ of index 1.
\par
Let $\xi $ be a local unit vector field normal to $M$ in $\Bbb CQ^m ,$  $A$ the shape operator of the
immersion $x,$ and let $\alpha = (1/n)\,\text{tr}\, A $ be the mean curvature of $M$ in $\Bbb CQ^m,$ so that the mean
curvature vector $H$ of the immersion equals $H = \alpha \xi .$
Further, let $\overline\nabla , \; \bar A, \overline D,$ denote respectively the Levi-Civita
connection, the shape operator, and the metric connection in the normal
bundle, related to $\cxf$ and the embedding $\Phi$. Let the same letters
without bar denote the respective objects for a submanifold $M$ and the immersion $x,$
whereas we use the same symbols with tilde to denote the corresponding objects related to the
composite immersion $\tilde x : = \Phi \circ x$ of $M$ into the (pseudo) Euclidean space
$H^{(1)}(m+1).$ As usual, we use $\s$ for the second fundamental form of
$\cxf$ in $E^N_{(K)}$ via $\Phi$ and $h$ for the second fundamental form of
a submanifold $M$ in $\cxf$. An orthonormal basis of the tangent space $T_pM$ at a general
point will be denoted by $\{e_i \}, i = 1, 2, {\ldots} , n.$ In general, indices $i, j$ will range from $1$ to $n$ and $\Gamma $ 
will denote the set of all (local) smooth sections of a bundle.
\par
We give first some important formulas
which will be repeatedly used throughout this paper.
For a general submanifold $M,\,$ local tangent fields $X, Y \in \Gamma (TM)$ and a
local normal field $\xi \in \Gamma (T^\perp M),$ the
formulas of Gauss and Weingarten are
$$ \overline \nabla_XY = \nabla_XY +
h (X, Y) ; \qquad\qquad \overline \nabla_X \xi = -A_\xi X + D_X\xi.  \tag3 $$

In particular, for a hypersurface of a complex space form $\Bbb CQ^m $
with unit normal vector $\xi $ and the corresponding shape operator $A,$ they
become
$$ \overline \nabla_XY =\nabla_XY + \langle AX, Y \rangle \xi
;\qquad \qquad \overline \nabla_X\xi = - AX \tag4 $$

Let $J$ be the K\"ahler almost complex structure of $\cxf$, and $U$ be the
distinguished tangent vector field $ U:= -J\xi .$ Define an endomorphism
$S$ of the tangent space and a normal bundle valued 1-form $F$ by
$$ SX = (JX)_T , \quad FX = (JX)_N  = \langle X, U \rangle \xi ,$$
i.e for $X \in \Gamma (TM), \;\; JX = SX + FX $ is the
decomposition of $JX$ into tangential and normal to submanifold parts.  Then the
following formulas are well known [3], [25]:

$$ SU = 0, \quad SX = JX - \langle X, U \rangle \xi , \qquad S^2X = -X
+ \langle X, U \rangle U \tag5
$$
$$ \nabla_XU = SAX, \qquad (\nabla_XS)Y = \langle Y, U \rangle AX -
\langle AX, Y \rangle U.  \tag6 $$
The curvature tensor of $\cxf (4c)$ is given by

$$ \overline R(X, Y)Z = c\, [ \la Y, Z\ra X - \la X, Z\ra Y +  \la JY , Z\ra JX -
\la JX , Z\ra JY - 2 \la JX , Y\ra JZ \, ], \tag7
$$
and the equations of Codazzi and Gauss for a hypersurface of $\Bbb CQ^m(4c)$ are
respectively given by
$$ (\nabla _XA)Y - (\nabla_YA)X = c\, [\langle X, U \rangle SY -
\langle Y, U \rangle SX - 2\langle SX, Y \rangle U], \tag8 $$

$$ \align R(X, Y)Z = \; &c\, [ \langle Y, Z \rangle X - \langle X, Z
\rangle Y + \langle SY, Z \rangle SX - \langle SX, Z \rangle SY -2
\langle SX, Y \rangle SZ ] \\ &+\langle AY, Z \rangle AX - \langle AX,
Z \rangle AY .\tag9 \endalign $$

The following formulas of A.  Ros for the shape operator of $\, \Phi
\,\, $ in the direction of $ \sigma (X,Y) $ are also well known, (see, for
example, [28], [29] and [18])
$$ \align \la \s (X,Y), \s (V, W) \ra = c\, [&2 \la X, Y\ra \la V, W \ra +
 \la X, V\ra \la Y, W \ra + \la X, W\ra \la Y, V \ra\\ &  +
\la JX, V\ra \la JY, W \ra + \la JX, W\ra \la JY, V \ra ], \tag 10 \endalign
$$

$$ \bar A_{\sigma (X,Y)}V = c\, [2 \langle X, Y\rangle V + \langle X,
V\rangle Y + \langle Y, V\rangle X + \langle JX, V\rangle JY + \langle
JY, V\rangle JX ]. \tag{11}
$$

One also verifies
$$\s (JX, JY) = \s (X, Y), \qquad \la \s (X, Y), \tilde x \ra = - \la X, Y
\ra , \qquad
\la \s (X, Y), \, I \ra = \, 0. \tag{12}
$$
The gradient of a smooth function $f$ is a vector field
$\nabla f := \sum_i (e_i f)e_i.$
The Hessian of $f$ is a symmetric tensor field defined by
$$  \text{Hess}\,_f(X, Y) = \la \nabla_X(\nabla f), Y \ra = XY f - (\nabla_XY)f,
$$
and the Laplacian acting on smooth functions is defined as $\lp f = -
\text{tr\, Hess}\,_f.$
The Laplace operator can be extended to act on a vector field $V$ along $\tilde x
(M)$ by
$$\lp V = \sum_i [\tilde\nabla_{\nabla_{e_i}e_i} V - \tilde\nabla_{e_i}
\tilde \nabla_{e_i} V].$$
The product formula for the Laplacian, which will be often used in the ensuing
computations, is
$$ \Delta (f\, g) = (\Delta f)\, g + f(\Delta g) - 2 \sum_i (e_if)( e_ig) ,
\tag{13} $$
for smooth functions $f, \,g \in C^\infty (M)\, ,$ and it can then be extended to hold for the scalar
product of vector valued functions, and thus also for product of matrices, in a natural way.  We shall use
the notation $f_k:= \text{tr}\, A^k, $ and in particular $f:= f_1 = \text{tr}\, A.$
For an endomorphism $B$ of the tangent space of $M$ we define
$\text{tr} (\nabla B): = \sum_{i=1}^n (\nabla_{e_i}B)e_i .$ We shall assume all manifolds to be smooth and
connected, but not necessarily compact.
\vskip 1truecm

\head 3. Iterated Laplacians of a real hypersurface \endhead

\vskip 0.6truecm
Recall that $$ \Delta \tilde x = - n \tilde H = - f\,\xi - \sum_{i=1}^n \sigma (e_i,
e_i ),  \tag{14} $$
where here, and in the following, we understand the Laplacian $\Delta $ of $M$ to be applied to
vector fields along $M$ (viewed as $E^N_{(K)}-$valued functions, i.e.  matrices) componentwise.

By the product formula above we have

$$ \Delta^2\tilde x := \Delta (\Delta\tilde x ) = - (\Delta f )\xi - f
(\Delta \xi ) + 2 \sigma (\nabla f , \xi ) - 2 A(\nabla f )- \sum_i
\Delta\, (\sigma (e_i, e_i )).  \tag{15} $$
Further,
$$ \align\Delta\xi = \sum_i &[\tilde \nabla _{\nabla_{e_i}e_i}\xi -
\tilde \nabla _{e_i}\tilde \nabla_{e_i} \xi ]\\ = \sum _i
&[-A(\nabla_{e_i}e_i ) + \sigma (\nabla_{e_i}e_i, \xi ) + \overline
{\nabla}_{e_i}(Ae_i) + \sigma (e_i, Ae_i)\\ & \;+ \bar A_{\sigma (e_i,
\xi )} e_i - \overline{D}_{e_i} (\sigma (e_i, \xi ))]. \endalign $$
Using (11), the parallelism of $\sigma ,$  and the fact that $\, \text{tr}\, (\nabla A) = \nabla (\text{tr}\, A) = \nabla f $
(by virtue of the Codazzi equation), we obtain

$$ \Delta\xi = \nabla f + [f_2 + c(n - 1)]\, \xi - f\, \sigma (\xi, \xi )
+ 2 \sum_i \sigma (e_i, Ae_i ).  \tag{16} $$

One further computes

$$ \align \sum_i \Delta ( \sigma (e_i, e_i)) = & -\, 4c\, JAU + 2c
(n+3)f\,\xi + 2c\, (n+2) \sum _i \sigma (e_i, e_i ) \\ & + 2 \sum_i
\sigma (Ae_i, Ae_i) - 2 \sigma (\xi, \nabla f) - 2(c + f_2 )\, \sigma
(\xi, \xi ).  \tag{17} \endalign $$

Combining formulas (15)-(17) we finally obtain

$$ \align \Delta^2 \tilde x = &- [\Delta f + f (f_ 2 + c( 3n + 5)) -
4c\langle AU, U \rangle ]\, \xi + 4c\, SAU \\ &- f\, \nabla f - 2
A(\nabla f) + ( 2c + 2 f_2 + f^2 )\, \sigma (\xi, \xi ) + 4\, \sigma
(\nabla f , \xi ) \\ &- 2 c(n+2) \sum_i \sigma (e_i, e_i ) - 2 f\,\sum_i
 \sigma (e_i, Ae_i ) -  2 \sum_i\sigma (Ae_i, Ae_i ) .
\tag{18} \endalign $$
Compare this with formula (2.15) of [36], formula (2.9) of [18], and formula (2.8) of [13].
\vskip 0.4truecm
Let us now find $\Delta^3\tilde x .$ The computation is long but
straightforward, so we just outline the main steps.  First we
shall compute the trace-Laplacian of the shape operator defined as the endomorphism
$\lp A := \sum_i[ \nabla_{\nabla_{e_i} e_i}A - \nabla_{e_i}(\nabla_{e_i}A)] $.
This computation is modeled on the computation of Nomizu and Smyth
[26] in the case of constant-mean-curvature-hypersurfaces of a real
space form.  However, here we do not assume the mean curvature to be
constant and we are dealing with complex space forms.
Let $K(X, Y) = \nabla_{\nabla_XY}A - \nabla_X(\nabla_YA).$
Then
$$ K(X, Y) = K(Y, X) + [A, R(X, Y)], \tag{19} $$
\noindent
where $\,\, R(X, Y) = \nabla_X \nabla _Y -
\nabla_Y \nabla_X - \nabla_{[X,Y]}\, $ is the curvature operator of
the hypersurface and the bracketed expression on the right hand side denotes
the commutator of the endomorphisms involved.  Clearly,
$\Delta A = \sum _i K(e_i, e_i).\, $ We compute

$$ \align \sum_i K(X, e_i) e_i &= \sum_i [( \nabla_{\nabla_Xe_i}A)e_i
- (\nabla_X (\nabla_{e_i}A)) e_i ]\\ &= \sum_{i,j} \omega_i^j(X) [(
\nabla_{e_j}A)e_i + ( \nabla_{e_i}A)e_j] - \sum_i \nabla_X
((\nabla_{e_i}A) e_i) \\ &= - \nabla_X (\nabla f) , \tag{20}\endalign $$
\noindent
since the connection 1-forms $\omega_i^j$ are antisymmetric
and the bracketed expression is symmetric in $i, j.$ Since all the
quantities involved are tensorial, to facilitate further computations
let us assume that $\nabla_{e_i}e_j = 0 \, $ for all $i,j = 1, ..., n $ and
moreover $\nabla_{e_i}X = 0 \, $ at a point where the computations are
being carried out.  By the Codazzi equation we have
$$ \align \sum_i K(e_i, X) e_i &= - \sum_i \nabla_{e_i}((\nabla_XA)e_i
) \\ &= - \sum_i \nabla_{e_i}[(\nabla_{e_i}A)X + c\, (\langle X, U
\rangle Se_i - \langle e_i, U \rangle SX - 2 \langle SX, e_i \rangle
U\, ) ]\\ &= \;\; \sum_i \, [K(e_i, e_i )X  + 2c\, \nabla_{e_i} ( \langle
SX, e_i \rangle U )\\
&\qquad\qquad\quad - c\, \nabla_{e_i}( \langle X, U
\rangle Se_i - \langle e_i, U \rangle SX )].  \endalign $$
Using (5) and (6) we get

$$ \align \sum_i \nabla_{e_i} ( \langle SX, e_i \rangle U ) &= \sum_i
(\,\langle (\nabla _{e_i}S)X, e_i \rangle U + \langle SX, e_i \rangle
\nabla_{e_i} U\, )\\ &= \;\;  f \langle X, U \rangle U -
\langle AX, U \rangle U
 + SASX, \endalign $$
and in a similar fashion

$$ \align \sum_i \, \nabla_{e_i}( \langle X, U \rangle Se_i - \langle
e_i, U \rangle SX ) &= \sum_i [\langle X, \nabla_{e_i}U \rangle Se_i + \langle
X, U \rangle (\nabla_{e_i}S) e_i ] \\ &\quad - \sum_i [\langle e_i,
\nabla_{e_i}U \rangle SX + \langle e_i, U \rangle (\nabla_{e_i}S) X] \\
&= \langle AX, U \rangle U - f \langle X , U \rangle U - SASX .
\endalign $$
Combining these steps we get
$$ \sum_i K(e_i, X) e_i = (\Delta A)X + 3c\, SASX -3c\, \langle AX -
f X ,\, U \rangle\, U .  \tag{21} $$
From (19) and (21) it follows
$$ \align (\Delta A)X =& - 3c\, SASX + 3c\,
\langle AX - f\, X ,\, U \rangle\, U \\ & \,\,\, + \sum_i K(X, e_i)e_i
+ \sum_i [A, R(e_i, X)] e_i .  \endalign $$
By the Gauss equation (9) we have
$$ \align \sum_i [A, R(e_i, X)]e_i \;\; &= \quad  \sum_i A(R(e_i, X)e_i) - \sum_i
R(e_i, X)(Ae_i) \\&= \quad  c f\, X + (f_2 - cn) AX - f A^2X + 3c\,AS^2X - 3c\, SASX .  \endalign $$
From (5) and the above we finally obtain the following extension of the
Nomizu-Smyth formula [26] to hypersurfaces of non-Euclidean complex space forms:
$$ \align (\Delta A)X = \quad &3c\, [\langle AU - f\, U , \, X \rangle\, U +
\langle U, X \rangle AU ] - 6 c\, SASX \\ &+ cf X + [f_2 - c(n+3)] AX -
f A^2 X - \nabla _X (\nabla f) .  \tag{22} \endalign $$
As in [26], we have $\lp (\text{tr}\, A^2) = 2\, \text{tr}[(\lp A) A] - 2
\sum_i \text{tr} (\nabla_{e_i}A)^2,$ and thus we obtain the following
Simons'-type formula for a hypersurface of $\Bbb CQ^m(4c):$
$$ \align \frac 1 2 \lp (\text{tr} \, A^2) = \; & 6c\, \la AU, AU\ra - 3cf
\la AU, U \ra - 6c\, \text{tr}\, (SA)^2 + cf^2\\& + [f_2 - c (n + 3)]\, f_2 - f
f_3 - \Vert\nabla A \Vert^2 - \sum_i \text{Hess}_f\, (Ae_i, e_i), \tag{23}
\endalign
$$
where $\Vert \nabla A \Vert^2 = \sum_{i, j} \la (\nabla_{e_j}A)e_i,(\nabla_{e_j}A)e_i  \ra . $
A similar, but rather long, calculation using (22) yields

$$ \align \Delta (JAU) &= [\text{tr}\, (\nabla_UA^2)]\, U - f SA^2U +
2 (c + f_2 ) SAU - 2 \sum_i J(\nabla_{e_i}A)(SAe_i) \\ &\quad - J \nabla_U (\n) - JAS (\n) + \la \n , AU\ra U
\\ &\quad - [ f \langle A^2U, U \rangle - 2 (c + f_2) \langle AU, U \rangle
+ 2c f]\, \xi \\ &\quad + 2 \sigma (A^2U, U)
- f \sigma (\xi , JAU) - 2 \sum_i  \sigma (e_i,\,
J\nabla_{e_i}(AU)).  \tag{24} \endalign $$
Note that by the Codazzi equation and formulas (5), (6) we have
$$ \nabla_{e_i}(AU) = ASA e_i + ( \nabla_UA) e_i - c \,
Se_i, \tag{25} $$
and also $\s (\xi , SX) = \s (U, X) - \la X, U\ra \s (\xi , \xi).$
Additional, somewhat involved, computations yield the following formulas:

$$ \align \lp ( \s (\xi , \xi)) &=  4c JAU + 2 (c + f_2)\, \s (\xi , \xi)
+ 2 \s (\n , \xi)\\ &\quad
+ 2c \sum_i \s(e_i , e_i) - 2 \sum_i \s (Ae_i, Ae_i),  \tag{26} \endalign
$$

$$\align \sum_i \Delta (\sigma (e_i, Ae_i)) = \;  &- 2c \sum_i J
(\nabla_{e_i}A)(Se_i) + 2cf JAU - 4c JA^2U + 8c \n \\ &+ 2 (cf^2 + 2c f_2 + n-1 )\, \xi
- 4\, \sigma (\xi, \text{tr}\, (\nabla A^2)) - 2 f_3 \, \sigma (\xi , \xi)\\
 &+ 2\, \sigma (\xi , A(\n)) + \sum_i \, [\sigma (e_i,
(\Delta A)e_i) + 2 \sigma (Ae_i, A^2e_i) ] \\ &+ 2c\, \sum_i [
f \sigma (e_i, e_i) +  \sigma (e_i, Ae_i) - \sigma (e_i, JASe_i)],
 \tag{27} \endalign $$
and using
$$ \sum_{i, j} \bar A_{\s ((\nabla_{e_j }A)e_i, A e_i)}e_j
= c\, \text{tr}\, (\nabla A^2) + c\, \n_2 - c \sum_j J (\nabla_{e_j} A^2)(Se_j),
$$
also
$$ \align \sum_i \Delta (\sigma (Ae_i, Ae_i)) &= - \, 2 \sum_{i,j} \sigma
((\nabla_{e_j}A)e_i, (\nabla_{e_j}A)e_i) + 2 \sum_i \sigma (Ae_i,
(\Delta A)e_i)  \\ & \quad + 2
\sum_i\, [cf_2 \,\sigma (e_i, e_i ) + c\, \sigma (e_i, A^2e_i) + \sigma (Ae_i, A^3e_i)]\\
&\quad - 2c\, \sum_i \sigma (e_i, JA^2Se_i) - 4\, \sigma (\xi, \text{tr}\, (\nabla A^3)) + 2\, \s (\xi , A^2(\n))\\
&\quad - 2 f_4 \, \sigma (\xi , \xi) + 4c\, \text{tr}\, (\nabla A^2) + 4c\, \nabla f_2  + 2c\, (f f_2 + 2f_3) \xi \\
&\quad + 2c f JA^2U - 4c\, JA^3U - 4c\, \sum_i J(\nabla_{e_i}A^2)(Se_i) .  \tag{28}\endalign $$
By a repeated use of the Codazzi equation we may deduce that
$$\text{tr}\, (\nabla A) = \nabla (\text{tr}\, A) = \nabla f, \tag{29}$$
$$ \text{tr}\, (\nabla A^2) = \frac{1}{2} \nabla f_2 + A(\nabla f) -
3c\, SAU, \tag{30}$$
$$ \text{tr}\, (\nabla A^3) = \frac{1}{3} \nabla f_3 + \frac{1}{2}
A(\nabla f_2) + A^2(\nabla f)- 3c\, SA^2U - 3c\, ASAU, \tag{31}$$
and in general, by induction,
$$ \text{tr}\, (\nabla A^k) = \sum_{r=1}^k \frac{1}{r} A^{k-r}(\nabla
f_r) - 3c\, \sum_{r=1}^{k-1}\, A^{r-1}SA^{k-r}U.  $$
Additionally, by using the Codazzi equation again one computes
$$ \sum_i (\nabla_{e_i}A)(Se_i) = - c (n-1) U, \tag{32}$$
and for a symmetric endomorphism $B$ one gets from (12)
$$ \sum_i \s (e_i , JB e_i)= 0 . \tag{33}
$$
Although we can compute $ \lp (\s (\n , \xi))$ in a similar fashion and obtain
a formula for $\lp^3\tilde x$ in general, we list this formula only in the special
case when the mean curvature is (locally) constant. Thus assuming $f = \, $ const, from
(18) and (22)-(33) we obtain
$$  \lp^3 \tilde x = (\lp^3 \tilde x)_T +  (\lp^3 \tilde x)_N,
$$
where the component tangent to $\cxf$ equals
$$\align  (\lp^3 \tilde x)_T  = \;\; & 8c \sum_i J [(\nabla_{e_i} A^2)(S e_i) -
(\nabla_{e_i}A)(SAe_i)] - 2 f A (\nabla f_2) - 4c\, \nabla f_2 \\
   & + 12 c\, \la \nabla f_2 , U\ra U + 8c\, SA^3U + 8 (2c f_2 + n + 7)\, SAU\\
   & + [8c \la A^3U, U\ra + 8 (2c f_2 + n + 4) \la AU, U\ra - f (\lp f_2) -
8c f_3\\
   & \quad\; - f (f_2^2 + 4c (n+4) f_2 + 4c f^2 + 7 n^2 + 30 n + 19) ]\, \xi ,
\tag{34} \endalign
$$
and the normal component is
$$\align (\lp^3 \tilde x)_N = & \, 4 \sum_{i, j} \s ( (\nabla_{e_j}A) e_i
, (\nabla_{e_j}A) e_i) + 6f \s (\xi , \nabla f_2) + 12 \s (\xi , A (\nabla
f_2))\\
   & + \frac 83 \,\s (\xi , \nabla f_3) - 16 c\, \s (\xi , ASAU) - 32 c f\, \s (AU
, U)\\
   & - 32 c\, \s (A^2U, U) - 12 c\ \s (AU, AU) + 16 c f \sum_i \s (e_i , SASe_i)\\
   & + 16 c \sum_i \s (e_i, SASAe_i) + 4c \sum_i \s (e_i , SA^2S e_i) \\
   & + \{ 28c\, \la A^2U, U\ra + 28 c f \la AU, U\ra + 2 (\lp f_2) +
f^2 [3 f_2 + c (3n + 13)]\\
   & \quad\;  + 4 f_2^2 + 4c (n+4) f_2 + 4 f f_3 + 4 f_4 + 4n + 20\} \, \s (\xi ,
\xi)\\
   & - 4\, (c f^2 + n^2 + 4n + 5) \sum_i \s (e_i, e_i) - 4f \, [f_2 + c
(n+3)] \, \sum_i \s (e_i, A e_i)\\
   & \quad -4 \, (c + 2 f_2) \, \sum_i \s (A e_i, A e_i) - 4 \, \sum_i \s (A^2
e_i, A^2 e_i). \tag{35}
\endalign
$$
Compare with formula (2.17) of [13] and a related expression in [17].
\vskip 1truecm
\head{4. Hopf hypersurfaces of 2-type have constant principal curvatures}
\endhead
\vskip 0.6truecm

In this section we study hypersurfaces of complex space forms of Chen-type 2 and
classify such hypersurfaces that are also assumed to be Hopf hypersurfaces i.e. for which the
structure vector field $U:= - J\xi$ is principal. We denote by $\Cal D^\perp$ the
1-dimensional distribution generated by $U$ and by $\Cal D$ the holomorphic distribution
which is the orthogonal complement of $\Cal D^\perp$ in $TM$ at each point. By way of notation, $V_{\mu }$ will denote
the eigenspace of the shape operator $A$ for an eigenvalue (principal curvature) $\mu$ and $\frak s (\Cal D),$ the spectrum of $A|_{\Cal D},$
the set of all eigenvalues of $A$ corresponding to eigenvectors belonging to $\Cal D$ at a given point.
\par
Let $M^n \subset \Bbb CQ^m, \, n = 2m - 1, $ be a 2-type hypersurface in
$H^{(1)}(m+1) ,$ i.e.
$\tilde x = \tilde x_0 + \tilde x_u + \tilde x_v $ where
$\tilde x_0 =\,$ const, $ \lp \tilde x_u  = \lambda_u \tilde x_u$ and
$\lp \tilde x_v = \lambda_v \tilde x_v, $ according to (1). Then
$$ \Delta ^2 \tilde x - (\lambda_u + \lambda_v )\Delta \tilde x +
 \lambda_u\lambda_v \tilde x = \lambda_u\lambda_v \tilde x_0 .  \tag{36}
$$
Let $L$ be the vector field in $H^{(1)}(m+1)$ along $M^n,$ represented by the left
hand side of the equation (36) and let $X$ be an arbitrary tangential vector
field of $M^n$.  Then

$$ \align 0 = \; & \la \tilde \nabla _XL, \tilde x \ra = X \la
L, \tilde x \ra - \la L, X \ra \\ = \; & X ( -2c + f^2 + 2cn (n+2) -
(\lambda_u + \lambda_v) n + \frac{c}{2} \lambda_u\lambda_v ) \\ & - 4c
\langle SAU, X \rangle + \langle f \nabla f + 2 A(\nabla f ), X
\rangle \\ = \; &\langle 2f \nabla f, X \rangle - 4c \langle SAU , X
\rangle + \langle f \nabla f + 2 A(\nabla f), X \rangle .  \endalign
$$
Therefore
$$ 2 A(\nabla f ) + 3 f \nabla f - 4c SAU = 0. \tag {37} $$
Similarly, by considering the $\sigma (\xi, \xi)-$component, in combination with
(37) we may obtain
$$ \align \nabla f_2 &+ A(\nabla f) - \frac{1}{2} f \nabla f - f
(\nabla_UA)U - 2 (\nabla_UA)(AU) \\ &+ f \langle \nabla f, U \rangle U
- 2 \langle \nabla f, U \rangle AU - \langle AU, U \ra \nabla f = 0,
\tag{38} \endalign $$
\noindent
and the other components are even more
complicated. Although it is possible to characterize 2-type hypersurfaces of
$\cxf$ by a set of equations involving the structure
vector field $U,$ the gradients of $f$ and $ f_2,\, \lp f, \,$
the shape operator, and various compositions of $S$ and $A,$ the equations involved are very complicated to enable the
classification of such hypersurfaces without any extra conditions.
At this point it seems beneficial to make some additional
assumptions on a hypersurface in order to make the situation more
tractable. The most facile assumption, which simplifies many terms, is
that $f: = \text{tr}\, A = \text{const},$  immediately leading, by way of
(37), to the conclusion that $M$ is a Hopf hypersurface, since
$SAU = 0$ is equivalent to $AU = \k\, U\, $ for some function
$\k .$ Moreover, it is known that in this case  $\k$ is (locally)
constant [22], [25].
Using this, one can show that the hypersurface is homogeneous and has at most 5 distinct
principal curvatures, all of which are constant. Using the complete list of such hypersurfaces
available in [32], [24], [20], [3], [4], one obtains a classification of constant-mean-curvature (CMC)
hypersurfaces whose  Chen-type is 2.  This has been already attempted by
Udagawa [36] for hypersurfaces of $\Bbb CP^m(4),$ and for hypersurfaces of
$\Bbb CH^m (-4)$ see below. Udagawa's classification in $\Bbb CP^m, $ however, is incomplete (see below).
On the other hand, instead of assuming the
mean curvature to be constant, it seems more challenging to make a weaker assumption
that $M$ is a 2-type Hopf hypersurface. In that case we have
$$ AU = \k\, U, \quad \k = \text{const}, \qquad\text{and} $$
$$ A(\nabla f) = - \frac{3f}{2} \nabla f. \tag {39} $$
So we do not get $f = \text{const}$ immediately, although that will eventually turn out to be the case.

Let $G$ be an open set defined by $G = \{ p \in M | \, f(p)\cdot (\nabla f) (p) \neq
0 \}.$ The Hopf property implies $\langle U, \nabla f\rangle = Uf = 0$ [25, p. 253] and thus
$\nabla f \in \Cal D = (\Bbb R U)^\perp .$ In addition, we have that $S(\nabla f)$
is also an eigenvector of $A,$ see [4], [22].
Then, since the integral curves of $U$ for a Hopf hypersurface are geodesics, (38) reduces to
$$ \nabla f_2 = (2f + \k)\,\nabla f. \tag{40} $$
Instead of showing more general formula (38), for our purposes it suffices to
prove (40). By using (3) and parallelism of $\s$ we have
$$ \align 0 &= \; \la \tilde\nabla_XL, \s (\xi, \xi )\ra\\
            &= \; X\la L, \s (\xi, \xi)\ra + \la L, \bar A_{\s (\xi, \xi)}X \ra
            - \la L, \overline{D}_X(\s (\xi, \xi))\ra \\
            &= \; X \,\la L, \s (\xi, \xi)\ra + 2c\, \la L, \, X + \la X, U\ra U \, \ra +
            2\la L, \s (AX, \xi)\ra.\tag{41} \endalign $$
Using $AU = \k U $ and $Uf = 0$ we obtain
$$ \la\lp^2 \tilde x , \, \s (\xi, \xi)\ra = 4c\, [f_2 - \k f - \k^2 - cn (n + 3)],
$$
$$ \la\lp^2 \tilde x , \, X + \la X, U\ra U \,\ra = - \la f \n + 2 A(\n ), X\ra ,\quad\;
\la\lp^2 \tilde x , \, \s (AX, \xi)\ra = 4c \la A(\n), X\ra .
$$
The metric products of $\lp\tilde x$ and $\tilde x$ with these quantities are either zero or give constants which disappear
after differentiation. Thus putting these together in (41) we obtain
$$  4c \la  \nabla f_2 - \k \n , X \ra - 2c  \la f \n + 2 A(\n ) , X\ra
+ 8c \la A(\n ), X\ra = 0,
$$
so that (40) follows from this and (39).
\par
We now show that $f = $\, const.
On $G$ let $e_1 := \n /|\n|$ be the unit vector of $\n .$ Then
$$ \align 0 &= \;\la\tilde\nabla_{\n}L, \s (e_1, e_1)\ra\\
            &= \; (\n) \la L, \s (e_1, e_1)\ra + \la L, \bar A_{\s (e_1, e_1)}\n \ra
                 - \la L, \overline{D}_{\n}\s (e_1, e_1)\ra \\
            &= \; (\n) \la L, \s (e_1, e_1)\ra + 4c\, \la L, \n \ra
               - 2 \,\la L, \s (\nabla_{\n}e_1, e_1)\ra
               + 3f |\n|\, \la L, \s (\xi, e_1)\ra.\endalign
$$
If $AX = \mu X$ for $X \in\Cal  D$ then, by the results of Maeda [22]
(for the projective case) and Berndt [4] (for the hyperbolic case), also $ A(SX) = \mu^*(SX),$ where
$\mu^*$ is uniquely determined by the condition
$$ (2\mu - \k)(2\mu^* - \k) = \k^2 + 4c, \quad \text{i.e.}\quad \mu^* =
\frac{\k\mu + 2c}{2\mu - \k}, \quad \text{and}\quad (\mu^*)^* = \mu . \tag{42}
$$
Since $\la U, \n \ra = 0, $ then $e_1, Se_1 \in \Cal D.$ When $X = \n ,$ we have $\mu = -3f/2 $ and
$\mu^* = \frac{3\k f - 4c}{2(3f + \k)}. $  We may assume that $3f + \k \neq 0,$ for otherwise we may work on an open subset
of $G$ where $f \neq - \k /3 $ and invoke continuity of $f.$ From here we compute using (11)
$$  \la L, \s (e_1, e_1)\ra = -c (5 f^2 + 4f\mu^* + 4\mu^{*2}) + \, \text{const},
$$
$$ \la L, \n\ra = 2f |\n|^2 , \quad
\la L, \s (\nabla_{\n}e_1, e_1)\ra = 0, \quad
\la L, \s (\xi, e_1)\ra = 4c |\n|.
$$
Substituting in the above equality we get
$$ 20 f|\n|^2 - (\n) (5 f^2 + 4f\mu^* + 4\mu^{*2}) = 0.
$$
Since
$$  (\n)(\mu^*) = \frac{3(\k^2 + 4c)}{2(3f + \k)^2}|\n|^2,
$$
we see that $f = \text{tr}\, A$ satisfies on $G$ a polynomial equation of degree 4
with constant coefficients, viz.
$$ 135 f^4 + 108 \k \,f^3 + 18 \k ^2 f^2 - 2 \k (5\k ^2 + 12c) f + 16 c\k^2 + 48 = 0.
$$
Consequently, $f$ is (locally) constant since $f$ is continuous and $M$ is assumed connected. From (40)
we also get $f_2 = \, $const.

Since $ f =\,$const, (18) reduces to
$$ \align \Delta^2 \tilde x = &\; [ 4c \k -  f (f_ 2 + c( 3n + 5))]\, \xi + ( 2c + 2 f_2 + f^2 )\, \sigma (\xi, \xi )\\
&\; - 2 c(n+2) \sum_i \sigma (e_i, e_i ) - 2 f\,\sum_i \sigma (Ae_i, e_i ) -  2 \sum_i\sigma (Ae_i, Ae_i ) .
\tag{43} \endalign $$
\noindent
 Differentiating (36) with respect to an arbitrary
tangent field  $X \in \Gamma (TM)$
we have
$$ \tilde\nabla_X (\Delta^2 \tilde x) - p \tilde\nabla_X (\Delta \tx ) + q X = 0, \tag 44
$$
where $p := \lambda_u + \lambda_v $ and $q := \lambda_u \lambda_v.$
Conversely, if (44) holds then $\tx$ satisfies the polynomial equation (36) in the Laplacian of the form $P(\lp )(\tx - \tx_0) = 0,$ with
$P(t) = t^2 - pt + q.$
According to a result of Chen and Petrovic [12] if such polynomial has simple real roots the submanifold is of 2-type (if not already of 1-type).
Let $V_\mu \subset \Cal D $ be an eigenspace of an eigenvalue $\mu \in \frak s (\Cal D)$ at each point and let $X \in V_\mu $ be a unit vector.
Taking the metric product of (44)
with $X$ and observing that $\langle \Delta \tx , X \rangle = \langle \Delta^2 \tx , X  \rangle = 0 \,$
for any tangent vector $X,$ we have
$$ \align 0 &= \langle \tilde\nabla_X (\Delta^2 \tx), X \rangle  - p \langle \tilde\nabla_X (\Delta \tx ), X \rangle  + q \\
&= X \la \Delta^2 \tx , X \ra - \la \Delta^2 \tx , \bar\nabla_X X + \s (X, X)\ra + p\, \la \Delta \tx , \bar \nabla_X X + \s (X, X)\ra + q \\
&= - \mu \la \Delta^2 \tx , \xi \ra - \la \Delta^2 \tx , \s (X, X)\ra + p \mu \la \Delta \tx , \xi \ra + p \la \Delta \tx , \s (X, X)\ra + q.
\endalign
$$
Using (10) and the above-mentioned results of Maeda and Berndt that $AX = \mu X $ implies $A (SX) = \mu^* (SX)$ where $\mu^* $ is
given by (42) we get from here
$$ \align 0 = q &+ [ f (f_2 + 3c (n + 3)) - pf - 4c \k ] \, \mu + 2c f^2 - 2 pc (n + 2)\\
&+ 4 (n + 1)(n + 3) + 4c \mu^2 + 4cf \mu^* + 4c {\mu^*}^2. \tag 45 \endalign
$$
Substituting the value of $\mu^*$ from (42) and clearing of denominators we get a fourth degree polynomial equation in $\mu$ with
constant coefficients (since $f, \, f_2, $ and $\, \k $ are all constant). We conclude that $A|_{\Cal D} $ has at most four eigenvalues,
i.e. the hypersurface has at most five distinct principal curvatures, all of them constant.
\vskip 0.4truecm
Hopf hypersurfaces of
$\Bbb CP^m(4)$ and $\Bbb CH^m(-4)$ for $m \geq 2 $ with constant principal curvatures are homogeneous and they are known. By a result of
Takagi [32] (see also [19], [20]) there are six types or  six classes of Hopf hypersurfaces with constant principal curvatures in $\Bbb CP^m(4),$
given as (possibly open portions of) the model hypersurfaces in the following list (the so-called Takagi's list):
\newline
($A_1$) \; A geodesic hypersphere of radius $r \in (0, \frac{\pi }{2});$
\newline
($A_2$) \; A tube of any radius $r \in (0, \frac{\pi }{2})$ around a canonically embedded (totally geodesic) $\Bbb CP^k$ for some
$k \in \{1, ..., m-2 \};$
\newline
($B$) \; A tube of any radius $r \in (0, \frac{\pi }{4})$ around a canonically embedded complex quadric
$Q^{m-1} = SO(m+1)/ SO(2) \times SO(m-1);$
\newline
($C$) \; A tube of radius $r \in (0, \frac{\pi }{4})$ around the Segre embedding of $\Bbb CP^1 \times \Bbb CP^k $ in $\Bbb CP^m, \; m = 2k+1;$
\newline
($D$) \; A tube of radius $r \in (0, \frac{\pi }{4})$ of dimension 17 in $\Bbb CP^9 $ around the Pl\"ucker embedding of the complex Grassmannian of 2-planes $G_2(\Bbb C^5);$
\newline
($E$) \; A tube of radius $r \in (0, \frac{\pi }{4})$ of dimension 29 in $\Bbb CP^{15}$ around the canonical embedding of the
Hermitian symmetric space $SO(10)/ U(5). $
\vskip 0.3truecm\noindent

We call these the standard examples or the model hypersurfaces in $\Bbb CP^m.$
To avoid confusion with the notion of Chen-type, these hypersurfaces will be referred to as being of class $A$ (with subclasses $A_1, \, A_2$),
$B, C, D, E,$ rather than being of type $A, B, C, D, E,$ as is customary in the literature.
For the example $B$ we note that a tube of radius $r \in (0, \frac{\pi }{4})$ around $Q^{m-1}$ in $\Bbb CP^m$ can be regarded
also as the tube of
radius $\frac{\pi }{4} - r$ around the canonically embedded (totally geodesic) $\Bbb RP^m$ in $\Bbb CP^m,$ which is the other focal 
submanifold of that hypersurface [5], [8].
\par
These model hypersurfaces have two, three, or five principal curvatures given by
$$ \k = 2 \cot (2r), \quad \text{and}\quad \mu_i = \cot (r + (i-1)\frac{\pi }{4}), \;\; i = 1, 2, 3, 4,
$$
where $r$ is the radius of the tube involved and $\k$ the principal curvature of $U.$
The table of principal curvatures and their multiplicities for these hypersurfaces is compiled by Takagi [32] and reads as follows
(see also [3], [25]):

\midinsert

\newdimen\tempdim
\newdimen\othick     \othick=.4pt
\newdimen\ithick     \ithick=.4pt
\newdimen\spacing    \spacing=9pt
\newdimen\abovehr    \abovehr=6pt
\newdimen\belowhr    \belowhr=8pt
\newdimen\nexttovr   \nexttovr=8pt

\def\rr{\hfil\down{\abovehr}&\omit\vrsp\vrule width\othick\cr
     \noalign{\hrule height\ithick}\up{\belowhr}&}
\def\up#1{\tempdim=#1\advance\tempdim by1ex
     \vrule height\tempdim width0pt depth0pt}
\def\down#1{\vrule height0pt depth#1 width0pt}
\def\large#1#2{\setbox0=\vtop{\hsize#1 \lineskiplimit=0pt \lineskip=1pt
     \baselineskip\spacing \advance\baselineskip by 3pt \noindent #2}
     \tempdim=\dp0\advance\tempdim by\abovehr\box0\down{\tempdim}}

\def\vrsp{\hskip\nexttovr\relax}
\def\toprule#1{\def\startrule{\hrule height#1\relax}}
\toprule{\othick}
\def\nstrut{\vrule height\spacing depth3.5pt width0pt}

\def\preamble#1{\def\startup{#1}}
\preamble{&##}
{\catcode`\!=\active
\gdef!{\hfil\vrule width0pt\vrsp\vrule width\ithick\relax\vrsp&}}

\def\table #1{\vbox\bgroup \setbox0=\hbox{#1}
\vbox\bgroup\offinterlineskip \catcode`\!=\active
\halign\bgroup##\vrule width\othick\vrsp&\span\startup\nstrut\cr
\noalign{\medskip}
\noalign{\startrule}\up{\belowhr}&}

\def\caption #1{\down{\abovehr}&\omit\vrsp\vrule width\othick\cr
\noalign{\hrule height\othick}\egroup\egroup \setbox1=\lastbox
\tempdim=\wd1 \hbox to\tempdim{\hfil \box0 \hfil} \box1 \smallskip
\hbox to\tempdim{\advance\tempdim by-20pt\hfil\vbox{\hsize\tempdim\noindent #1}\hfil}\egroup}

$$\table{\bf Table 1}
  ! \hfil $2\cot(2r)$ ! \hfil $\cot r$ ! $\cot (r + \frac{\pi}{4})$ ! $\cot (r + \frac{\pi}{2})$ ! $\cot (r + \frac{3\pi}{4})$  \rr
$\hfil A_1$ ! \hfil 1 ! \hfil $2(m-1)$ ! \hfil --- ! \hfil --- ! \hfil --- \rr
$\hfil A_2$ ! \hfil 1 ! \hfil $2(m - k - 1)$ ! \hfil --- ! \hfil $2k$ ! \hfil --- \rr
$\hfil B$ ! \hfil 1 ! \hfil --- ! \hfil $m-1$ ! \hfil --- ! \hfil $m-1$ \rr
$\hfil C$ ! \hfil 1 ! \hfil $m-3$ ! \hfil 2 ! \hfil $m-3$ ! \hfil 2 \rr
$\hfil D$ ! \hfil 1 ! \hfil 4 ! \hfil 4 ! \hfil 4 ! \hfil 4 \rr
$\hfil E$ ! \hfil 1 ! \hfil 8 ! \hfil 6 ! \hfil 8 ! \hfil 6\hfil
\caption{\sl Principal curvatures of the standard examples in $\Bbb CP^m$ and their multiplicities }
$$

\endinsert

It is known that the almost complex structure $J$ leaves eigenspaces $V_{\mu_1}$ and $V_{\mu_3}$ invariant and interchanges
eigenspaces $V_{\mu_2}$ and $V_{\mu_4}.$

In the complex hyperbolic space the number of principal curvatures is two or three. The list (the so-called Montiel's list after [24], completed  by Berndt [3], [4], see also [25]) of Hopf
hypersurfaces with constant principal curvatures in $\Bbb CH^m(-4)$ consists of (open portions of) the following:
\vskip 0.3truecm\noindent
($A_0$) \; A horosphere in $\Bbb CH^m ;$
\newline
($A^{\prime}_1$) \; A geodesic hypersphere of any radius $r \in \Bbb R_+;$
\newline
($A^{\prime\prime}_1$) \; A tube of any radius $r \in \Bbb R_+$ over a totally geodesic complex hyperbolic hyperplane $\Bbb CH^{m-1};$
\newline
($A_2$) \; A tube of any radius $r \in \Bbb R_+$ about the canonically embedded $\Bbb CH^k$ in $\Bbb CH^m$ for $ k = 1, ..., m-2$;
\newline
($B$) A tube of any radius $ r \in \Bbb R_+ $ about the canonically embedded (totally geodesic, totally real) $\Bbb RH^m$ in $\Bbb CH^m.$
\vskip 0.3truecm
We note  that a canonically embedded $\Bbb RH^m \subset \Bbb CH^m$ is of 1-type in $H^1(m+1), $ [15].
In a recent work  Berndt and D\'iaz-Ramos [6], [7] classified hypersurfaces of $\Bbb CH^m$ with three
constant principal curvatures, without assuming them to be Hopf.
\par
The table of principal curvatures $\k, \, \mu , \; \nu$ and their multiplicities $m_{\varkappa}, m_{\mu}, m_{\nu} $ 
is as follows [3], [4], [25]:
\vskip 0.2truecm
\midinsert

\newdimen\tempdim
\newdimen\othick     \othick=.4pt
\newdimen\ithick     \ithick=.4pt
\newdimen\spacing    \spacing=9pt
\newdimen\abovehr    \abovehr=6pt
\newdimen\belowhr    \belowhr=8pt
\newdimen\nexttovr   \nexttovr=8pt

\def\rr{\hfil\down{\abovehr}&\omit\vrsp\vrule width\othick\cr
     \noalign{\hrule height\ithick}\up{\belowhr}&}
\def\up#1{\tempdim=#1\advance\tempdim by1ex
     \vrule height\tempdim width0pt depth0pt}
\def\down#1{\vrule height0pt depth#1 width0pt}
\def\large#1#2{\setbox0=\vtop{\hsize#1 \lineskiplimit=0pt \lineskip=1pt
     \baselineskip\spacing \advance\baselineskip by 3pt \noindent #2}
     \tempdim=\dp0\advance\tempdim by\abovehr\box0\down{\tempdim}}

\def\vrsp{\hskip\nexttovr\relax}
\def\toprule#1{\def\startrule{\hrule height#1\relax}}
\toprule{\othick}
\def\nstrut{\vrule height\spacing depth3.5pt width0pt}

\def\preamble#1{\def\startup{#1}}
\preamble{&##}
{\catcode`\!=\active
\gdef!{\hfil\vrule width0pt\vrsp\vrule width\ithick\relax\vrsp&}}

\def\table #1{\vbox\bgroup \setbox0=\hbox{#1}
\vbox\bgroup\offinterlineskip \catcode`\!=\active
\halign\bgroup##\vrule width\othick\vrsp&\span\startup\nstrut\cr
\noalign{\medskip}
\noalign{\startrule}\up{\belowhr}&}

\def\caption #1{\down{\abovehr}&\omit\vrsp\vrule width\othick\cr
\noalign{\hrule height\othick}\egroup\egroup \setbox1=\lastbox
\tempdim=\wd1 \hbox to\tempdim{\hfil \box0 \hfil} \box1 \smallskip
\hbox to\tempdim{\advance\tempdim by-20pt\hfil\vbox{\hsize\tempdim\noindent #1}\hfil}\egroup}

$$\table{\bf Table 2}
  ! \hfil $\varkappa$ ! \hfil $\mu$ ! \hfil $\nu$ ! \hfil $m_{\varkappa}$ ! \hfil $m_{\mu }$ ! \hfil $m_{\nu }$ \rr
\hfil $A_0$ ! \hfil 2 ! \hfil --- ! \hfil 1 ! \hfil 1 ! \hfil --- ! \hfil $2m-2$ \rr
\hfil $A^{\prime}_1$ ! \hfil $2\coth (2r)$ ! \hfil $\coth r $ ! \hfil --- ! \hfil 1 ! \hfil $2(m-1)$ ! \hfil --- \rr
\hfil $A^{\prime\prime}_1$ ! \hfil $2\coth (2r)$ ! \hfil --- ! \hfil $\tanh r$ ! \hfil 1 ! \hfil --- ! \hfil $2(m-1)$ \rr
\hfil $A_2$ ! \hfil $2\coth (2r)$ ! \hfil $\coth r $ ! \hfil $\tanh r$ ! \hfil 1 ! \hfil $2(m- k-1)$ ! \hfil $2k$ \rr
\hfil $B$  ! \hfil $2\tanh (2r)$ ! \hfil $\coth r $ ! \hfil $\tanh r$ ! \hfil 1 ! \hfil $m-1$ ! \hfil $m-1$ \hfil
\caption{\sl Principal curvatures of the standard examples in $\Bbb CH^m$ and their multiplicities}
$$

\endinsert

It is known that the eigenspaces $V_{\mu}$ and $V_{\nu}$ are interchanged by the action of  $J$ for a class-$B$ hypersurface
and they are $J-$invariant (holomorphic) for any of the class-$A$ hypersurfaces. A hypersurface of class $A_2$ has three principal 
curvatures and so does a hypersurface of class $B$, except in one case, namely when the radius of the tube is 
$r = \frac{1}{2} \ln (2 + \sqrt 3 )$ and then $\mu = \varkappa = \sqrt 3. $
\par
In both settings, $\k$ is the principal curvature corresponding to $U:= - J\xi .$ These classifications enable us to prove our results
for 2-type Hopf hypersurfaces. In the subsequent investigation of 2-type Hopf hypersurfaces of $\Bbb CQ^m(4c)$ we  may assume
we are dealing with Hopf hypersurfaces with constant principal curvatures and therefore with one from the Takagi's list in the projective space
or one from the Montiel's list in the hyperbolic space. 
\vskip 0.8truecm
\head{5. The classification of 2-type Hopf hypersurfaces of $\Bbb CQ^m(4c)$}
\endhead
\vskip 0.3truecm

We begin by analyzing various components of equation (44).
Let $X \in \Gamma (TM).$ Then using the Gauss and Weingarten formulas (3) and the fact that  $\sigma$ is parallel, we get from (14)
$$ \tilde\nabla_X (\lp\tx ) = 2c (n + 2) X + f AX - 2 c \la X, U \ra U - f \s (X, \xi ) - 2 \s (AX, \xi ) \tag 46
$$
and from (43)
$$ \align \tilde\nabla_X (\lp^2 \tx ) = \; &[2c f^2 + 4 (n + 1)(n + 3)] X - [4c \k - f (f_2 + 3c (n + 3))]\, AX\\
& + 4c A^2 X - 2c\, [2 f_2 + f^2 + 2c (n + 3)]\la X, U\ra\, U - 4cf J ASX\\
& - 4c J A^2 SX  + [4c \k - f (f_2 + c (3n + 5))]\, \s (X, \xi )\\
& - 2\, [2f_2 + f^2 + 2c (n + 3)]\, \s (AX , \xi) - 4 f \s (A^2X, \xi )\\
& - 4 \s (A^3X, \xi ) - 2 f \sum_i \s ((\nabla_XA)e_i , e_i) -  4 \sum_i \s ((\nabla_XA)e_i , Ae_i).
\endalign
$$
Therefore, separating the part of equation (44) that is tangent to $\Bbb CQ^m$ we get
$$ \align 0 = & -2c\, [2f_2 + f^2 + 2c (n + 3) - p ] \la X, U\ra U - 4cf JA SX \\
& - 4c JA^2 SX + [4 (n + 1)( n + 3) + 2c f^2 - 2pc (n + 2) + q] \, X \\
& - [4c \k - f (f_2 + 3c (n + 3) - p)]\, AX  + 4c\, A^2X \tag 47
\endalign
$$
and the part normal to $\Bbb CQ^m$ yields
$$ \align & [4c\k - f (f_2  + c (3n + 5) - p)]\, \s (X, \xi) - 4 \s (A^3X, \xi )\\
& - 4f \s (A^2X, \xi ) - 2\, [2f_2 + f^2 + 2c (n + 3) - p] \, \s (AX, \xi ) \\
&- 2 f \sum_i \s ((\nabla_XA) e_i , e_i) - 4 \sum_i \s ((\nabla_XA)e_i , Ae_i) = 0. \tag 48 \endalign
$$
These expressions are linear in $X.$ Further separation of parts relative to the splitting
$\Cal D\oplus \Bbb R U  \oplus \Bbb R \xi $ of the tangent space of $\Bbb CQ^m$ yields the following

\proclaim{Lemma 1} Let $M^n$ be a Hopf hypersurface (not necessarily compact) of $\Bbb CQ^m(4c)$ $ \, (m \geq 2, \; n = 2m-1).$ If $M$
is of 2-type via $\tx$ satisfying 2-type condition (44) then $M$ has at most five distinct principal curvatures, all of which are constant,
and  the following relations hold:
\roster
\item "{($E_1$)}" \, $[2c (n + 1) + \k f]\, p =  q + \k f [f_2 + 3c (n + 3)] - 4c f_2 + 4n (n + 3);$
\item "{($E_2$)}" \, $$ \align  [2c (n + 2) + \mu f]\, p = \;\;& q + 2cf^2 + 4cf \mu^* + 4c {\mu^*}^2 + 4c \mu^2\\
  & + [f (f_2 + 3 c (n + 3)) - 4c \k ]\, \mu + 4 (n + 1)(n + 3),  \endalign $$ for any principal curvature $\mu \in \frak s (\Cal D) ;$
\item "{($E_3$)}" \, $$\align (f + 2 \mu )\, p = \;& - 4c\k + 4 \mu (f \mu + \mu^2 + f \mu^* + {\mu^*}^2)\\
  & + 2 \mu \, [2 f_2 + f^2 + 2c (n + 3) - 2f \k - 2 \k^2] + f [f_2 + c (3n + 5)],\endalign $$ for any $\mu \in \frak s (\Cal D) ;$
\item "{($E_4$)}" \, $$ \align &f \la (\nabla_XA)Y , Z \ra + f \la (\nabla_XA)(SY) , SZ \ra\\
& +  \la (\nabla_XA^2)Y , Z \ra + \la (\nabla_XA^2)(SY) , SZ \ra = 0,
\endalign $$
for every $X, Y, Z \in \Gamma (\Cal D).$
\endroster
Conversely, if $(E_1) - (E_4)$ hold for a Hopf hypersurface with constant principal curvatures, where $p$ and $q$ are
constants and $\mu \in \frak s (\Cal D)$ is an arbitrary principal curvature on $\Cal D,$ then the formula (44) holds and the submanifold
is of type  $\leq 2$ if the corresponding monic polynomial
$P(t) = t^2 - p t + q $ has two distinct real roots.

\endproclaim
\demo{Proof}
From the above  discussion it follows that $M$ has constant principal curvatures.
($E_1$) follows from (47) when $X = U.$ ($E_2$) is the formula (45), and it follows from (47)  when $X \in V_\mu \subset \Cal D$ is
chosen to be a principal direction  of a  principal curvature $\mu \in \frak s (\Cal D).$  Note that (47) is linear in $X,$ so it 
suffices to consider $X$ to be one of the principal directions.

 ($E_3$) and ($E_4$) follow from the normal part (48). Recall that the normal space of
$\Bbb CQ^m $ in $H^{(1)}(m + 1)$ is spanned by $\tx $ and the values of $\s$ on various pairs of tangent vectors to
$\Bbb CQ^m, \, $ namely $\s (\xi , \xi), \, \s (\xi , X), \, \s (X, Y)$ for $X, Y \in \Gamma (\Cal D), $ see e.g. [15]. Note that by (12)
$\s (\xi , U) = 0, \; \s (U, U) = \s (\xi , \xi)$ and $\s (U, X) = \s (\xi , JX)$ for $X \in \Cal D.$ By (12) and the constancy of
$f$ and $f_2,$ from (48) we conclude that the equation (44) has no $\tx-$component.
\par
Let $L$ denote the left-hand side of (48). Then $\bar A_L = 0.$ Conversely, if $\bar A_L = 0,$ for some
$L \in T^\perp \Bbb CQ^m$ then $L = k I$ is a multiple of the identity, but since $L$ is a linear combination of terms of the form
$\s (V, W)$ then by (12) $L = 0. $ Consider first $\bar A_L \xi = 0.$ The condition $\la \bar A_L\xi , \xi \ra = 0 $ gives no information
since by (10) $\la L, \s (\xi , \xi)\ra = 0$ is trivially satisfied and the same holds for  $\la \bar A_L\xi , U \ra = 0. $ 
Now take $Y \in \Cal D$ and consider  $\la \bar A_L\xi , Y \ra = \la L, \s (\xi , Y)\ra = 0. $ Using
$$\bar A_{\s (X, \xi)}\xi = c (X - \la X, U\ra U),\qquad \sum_i \bar A_{\s ((\nabla_XA)e_i , e_i)}\xi = 2c \k\, JSAX - 2 c\, JASAX,
$$
$$ \text{and} \quad  \sum_i \bar A_{\s ((\nabla_XA)e_i , Ae_i)}\xi = c \k ^2\, JSAX -  c\, JA^2SAX, $$
from (48) it follows
$$ \align &4 \la SA^2SAX, Y \ra + 4f \la SASAX , Y \ra - 4 \la A^3X, Y \ra - 4f \la A^2X, Y\ra\\
& -2 \, [2 f_2 + f^2 + 2c (n + 3) - p - 2f \k - 2 \k^2] \la AX, Y\ra\\
& + [4c \k + f p - f (f_2 + c (3n + 5))] \la X, Y\ra = 0.
\tag 49 \endalign $$
Since  $A\Cal D \subset \Cal D, \; S \Cal D = \Cal D,$ and  the expression
is linear in $X, Y \in \Cal D$ we can drop $Y$ and take $X \in V_\mu \subset \Cal D \,$ to get ($E_3$). Considering $\bar A_LU = 0$
gives no additional information beyond ($E_3$) by virtue of  $\la \bar A_LU, Y\ra = \la \bar A_L\xi, JY\ra ,$ returning it to the case above.
Next we exploit the condition $\bar A_L Y = 0 $ for $Y \in \Cal D.$ By (11) we have
$$\sum_i \bar A_{\s ((\nabla_XA)e_i ,\, e_i)}Y = 2c \, (\nabla_XA)Y - 2c\, J (\nabla_XA)(SY)$$
$$\sum_i \bar A_{\s ((\nabla_XA)e_i ,\, Ae_i)}Y = c\, (\nabla_XA^2)Y - c\, J (\nabla_XA^2)(SY).$$
In particular, when $X = U$ by the Codazzi equation and formula (6) we have
$$\sum_i \bar A_{\s ((\nabla_UA)e_i , e_i)}Y = 4SY + 2c\, [\k SAY + \k ASY - ASAY + (SA)^2SY]
$$
and a similar, somewhat longer, expression is obtained for $\sum_i \bar A_{\s ((\nabla_UA)e_i , Ae_i)}Y.\, $  Then taking  $Y \in V_\mu$ and
using $2 \mu \mu^* = 2c + \k (\mu + \mu^* ),$ by way of (42), we see that $\bar A_LY = 0$ reduces to a trivial identity when $X = U$.
Thus consider as the
last condition to check $\la \bar A_LY , Z\ra = 0.$ Choosing $Z = \xi $ or $Z = U$ gives back ($E_3$) and when $X, Y, Z \in \Cal D$
from $\la L, \s (Y, Z)\ra = 0$ we get ($E_4$).
Conversely, since we considered all possible components, the conditions ($E_1$) - ($E_4$) are equivalent to (47) and (48) by linearity
and thus we get (44), from which it follows that a hypersurface is of type $\leq 2,$ provided that the corresponding polynomial has two
distinct real roots.
\enddemo
\vskip 0.1truecm

Note that by a result of Niebergall and Ryan [25, p. 264] any of the class$-A$ hypersurfaces in $\Bbb CQ^m $ from either list  is
characterized by
$$(\nabla_XA)Y = - c\, [\la SX, Y \ra U + \la U, Y \ra SX],$$
so that the condition ($E_4$) is trivially
satisfied for those hypersurfaces. Further, by eliminating $q$ from ($E_1$) and ($E_2$) we get
$$ \align [2c + f (\mu - \k)]\, p = \; & 4c (\mu^2 + {\mu^*}^2) + 4cf \mu^* - 4c \k \mu + 4 (n + 3)\\
& + 4c f_2 + 2c f^2 + f (\mu - \k) [f_2 + 3c (n + 3)], \tag 50
\endalign $$
and if $p$ can be uniquely determined from this condition (regardless of the choice of $\mu$ and consistent with ($E_3$) ) then $q$ is
uniquely determined from ($E_1$).

\par
We now examine which of the Hopf hypersurfaces with constant principal curvatures are of 2-type. This has been already considered
by Udagawa for hypersurfaces of $\Bbb CP^m$ [36]. Although our argument is different from Udagawa's and relies on the analysis
of the conditions ($E_1$) - ($E_4$), rather than on the matrix representation of the immersion in $H^{(1)}(m+1)$, it partly overlaps
Udagawa's investigation and reaches the same classification for 2-type CMC real hypersurfaces in $\Bbb CP^m$ of class $A.$ However,
Udagawa's paper contains errors regarding mass-symmetric hypersurfaces and in particular hypersurfaces of class $B,$ 
as a result of which the three theorems in that work contain inaccuracies and incomplete classifications. Moreover, our more detailed analysis clearly exhibits the manner
of 2-type decompositions involved. Also, the benefit of our uniform approach is that it produces results for hypersurfaces of $\Bbb CH^m$ at the same time,
the case which is not treated in earlier papers, and the same technique will be used to study 3-type submanifolds.
\par
First we note that a horosphere in $\Bbb CH^m$ is not of any finite type since, as shown in [18], it satisfies
$\lp ^2 \tx = \, \text {const}\, \neq 0$ and therefore cannot satisfy equation (1), for otherwise equation (44) would hold for 
some constants $p$ and $q,$ which would force $p$ and $q,$ and thus also $\lp^2\tx ,$ to be zero or
$p \tilde\nabla_X (\lp \tx )$ to be a multiple of  $X,$ contradicting (46).
For hypersurfaces of class $A_1$ (geodesic spheres, equidistant hypersurfaces) we have

\proclaim{Lemma 2} \, $(i)$ \, A geodesic hypersphere in $\Bbb CP^m(4)$ of any radius $r \in (0, \pi /2), $  
$ \; r \neq \cot^{-1}\sqrt{1/ (2m + 1)}$ is of 2-type in $H(m+1).$
A geodesic hypersphere in $\Bbb CH^m(-4)$ of arbitrary radius $r > 0$ is of 2-type in $H^1(m + 1)$ via $\tx $ and the same holds true for a 
tube of an arbitrary radius $r > 0 $ about a totally geodesic complex hyperbolic hyperplane $\Bbb CH^{m-1}(-4) \subset \Bbb CH^m(-4).$  
These statements are also valid for any open portion of the respective submanifolds.
\vskip 0.2truecm\noindent
$(ii)\, $ The only complete mass-symmetric hypersurfaces of class $A_1$ are geodesic hyperspheres of radius $r = \cot^{-1}\sqrt{1/m}$
in $\Bbb CP^m(4).$
\endproclaim
\vskip 0.3truecm\noindent
\demo{Proof} \, $(i) \; $For a geodesic sphere (class $A_1$ in $\Bbb CP^m$ and $A^{\prime}_1$ in $\Bbb CH^m$ ) define
$$ \cot_c(r) = \cases \cot r, \; \; \;\text{when} \; c = 1 \; \text{(projective case)}\\
          \coth r, \; \text{when} \; c = - 1 \; \text{(hyperbolic case)}\endcases
$$
and let $\mu = \cot_c (r)$  be the principal curvature of multiplicity $2(m-1) = n-1$ and $\k = 2 \cot_c (2r)$ the principal curvature
(of $U$) of multiplicity 1, whereas $\mu = \tanh r, \; \k = 2 \coth (2r)$ for a tube about a complex hyperbolic hyperplane $\Bbb CH^{m-1}(-4)$
 of class $A^{\prime\prime}_1.$ Then
$$ \mu^* = \mu , \quad \k = \mu - \frac{c}{\mu}, \quad f = n \mu - \frac{c}{\mu}, \quad f_2 = n \mu^2 + \frac{1}{\mu^2} - 2c .\tag 51
$$

From (50) we get
$$ [(n + 2)c - \mu^{-2}]\, p = c (3n+2) (n + 2) \mu^2 + (3n^2 + 6n + 4) - \frac{(2n + 1)c}{\mu^2} - \frac{1}{\mu^4}. \tag 52
$$
We may assume that $(n+2)c \neq 1/\mu^2 ,$ certainly true when $c = -1,$ and when $c = 1$ the equality would lead to
$\mu = \sqrt{1/(n+2)}$ i.e. to $r = \cot^{-1}\sqrt{1/(2m+1)}.$ However, the geodesic hypersphere of this radius in $\Bbb CP^m(4)$
is of 1-type (see e. g. [23], [15]). Thus dividing (52) by $(n+2)c - \mu^{-2}$ we get
$$ p = (3n + 2) \mu^2 + 3c (n + 1) + \frac{1}{\mu^2} = (\mu^2 + c) (3n + 2 + \frac{c}{\mu^2}). \tag 53
$$
Then from ($E_1$) we find
$$ q = 2(n + 1)\left[ n \mu^4 + c (2n+1)\mu^2 + \frac{c}{\mu^2} + (n + 2)\right] .\tag 54
$$
Solving ($E_3$) for $p$  gives the same value as in (53), so the conditions $(E_1)-(E_3)$ are consistent and satisfied by the above
values of $p$ and $q,$ the condition ($E_4$) being trivially satisfied. According to Lemma 1, the equation (44) then holds, hence also (36).
Moreover the polynomial $P(\lambda ) = \lambda^2 - p \lambda + q $ has two distinct real roots $\lambda_u = 2 (n + 1)(\mu^2 + c)$ and
$\lambda_v = \frac{1}{\mu^2}(\mu^2 + c)(n\mu^2 + c),$ which are the two eigenvalues of the Laplacian from the 2-type decomposition of $\tx .$
\vskip 0.3truecm\noindent
$(ii)\;$ Let $\Cal D$ be the holomorphic distribution in $TM$ as before and choose an orthonormal basis $\{e_i \}$ of the tangent space so that 
$e_{n} = U$ and $e_i \in \Cal D$ for $i = 1, 2, ..., n-1.$ To see which hypersurfaces of class $A_1$ are mass-symmetric first we find 
from (14), (18), and (43)
$$ \Delta \tx = -(n \mu - \frac{c}{\mu}) \xi - \s (\xi , \xi) - \sum_{e_i \in \Cal D} \s (e_i, e_i), \quad \text{and}
$$
$$ \align \Delta^2 \tx =& - [n^2 \mu^3 + c (3n^2 + 2n - 4) \mu - \frac{2n-1}{\mu } - \frac{c}{\mu^3}] \, \xi\\ 
& + [(n^2 - 2) \mu^2 - 2cn - \frac{1}{\mu^2}] \, \s (\xi , \xi) - 2 (n+1)(\mu^2 + c) \sum_{e_i \in \Cal D} \s (e_i, e_i)\endalign
$$
and then compute, using $\lambda_u - \lambda_v = (1/\mu^2)(\mu^2 + c)[(n+2)\mu^2 - c],$ that
$$ \align \tx_u &= \frac{1}{\lambda_u (\lambda_u - \lambda_v)} (\Delta^2 \tx - \lambda_v \Delta\tx) \\
&= \frac{m-1}{4m (\mu^2 + c)^2} \{ - 4c \mu \xi + 2 \mu^2 \s (\xi , \xi)  
- \frac{\mu^2 + c}{m-1} \sum_{e_i \in \Cal D} \s (e_i, e_i)\} \tag 55\endalign
$$
$$  \tx_v = \frac{1}{\lambda_v (\lambda_v - \lambda_u)} (\Delta^2 \tx - \lambda_u \Delta\tx)
= - \frac{\mu}{(\mu^2 + c)^2}[(\mu^2 - c) \xi + \mu \s (\xi , \xi)]. \tag 56
$$
From Lemma 1 of [15] we have
$$ \tx = \frac{I}{m+1} - \frac{c}{2(m+1)} \s (\xi , \xi) - \frac{c}{4 (m+1)} \sum_{e_i\in \Cal D} \s (e_i, e_i).\tag{57}
$$
Now the center of mass can be found as $\tx_0  = \tx - \tx_u - \tx_v $ to yield
$$ \tx_0 = \frac{I}{m+1} + \frac{m\mu^2 -  c}{m(\mu^2 + c)^2} [\mu \xi + \frac{1}{2} \s (\xi , \xi) + (\mu^2 + c) (\tx - \frac{I}{m+1})] 
\tag{58}
$$
We observe that the same formula applies also for the center of mass of a 1-type hypersphere in $\Bbb CP(4)$ for an appropriate 
value of $\mu .$ Note that by our definition of mass-symmetry, any null 2-type hypersurface is per force mass-symmetric,
since the constant part $\tx_0$ can be manipulated and changed to be equal to $I/(m+1),$ and the existing constant $\tx_0$ moved to be a part 
of the $0-$eigenfunction. However, in our case, for $A_1$ hypersurface both $\lambda_u$ and $ \lambda_v$ as given above are nonzero,
because (in the hyperbolic case) $\mu = \coth r > 1.$
Since the $\xi$- component of the right hand side of (58) is the only part tangent to $\Bbb CQ^m$ and $\mu = \cot_cr \neq 0,$ a class-$A_1$
hypersurface is mass-symmetric, i.e. $\tx_0 = I/(m+1)$ if and only if $m \mu^2 = c.$ This is possible only when $c = 1$ and 
$\mu = \cot r = \sqrt{1/m}.$ Thus a geodesic hypersphere in $\Bbb CP^m (4)$ of radius $r = \cot^{-1}\sqrt{1/m}$ is the only complete
mass-symmetric hypersurface of class $A_1.$ We observe that this hypersphere with the given radius does satisfy the equation (3.14)
of [36], but it is completely overlooked in that paper. 
\enddemo

\proclaim{Lemma 3} $(i)\; $ There are no 2-type hypersurfaces in $\Bbb CH^m(-4)$ of class $A_2,$ i.e. no 2-type tubes about canonically 
embedded $\Bbb CH^k \subset \Bbb CH^m, \; 1 \leq k \leq m-2.$ A hypersurface of class $A_2$ in $\Bbb CP^m(4)$ is of 2-type if and only 
if it is an open portion of either (a) the tube of radius $r = \cot^{-1}\sqrt{\frac {k+1}{m-k}}$ or (b) the tube of radius 
$r = \cot^{-1}\sqrt{\frac{2k+1}{2(m-k)+ 1}},$ about a canonically embedded, totally geodesic  $\, \Bbb CP^k(4) \subset \Bbb CP^m(4),$ 
for any $k = 1, 2, ..., m-2.$
\vskip 0.3truecm\noindent
$(ii)\; $ The only complete mass-symmetric 2-type hypersurfaces of class $A_2$ are those in the first series of tubes (a) above.
\endproclaim

\demo{Proof} $(i)\; $Let $\mu_1 = \cot r, \; \mu_3 = \cot(r + \frac{\pi}{2}) = - \frac{1}{\mu_1}$ for model hypersurface of class 
$A_2$ in $\Bbb CP^m$
and $\mu_1 = \mu =  \coth r, \; \mu_3 = \nu = \tanh r =  \frac{1}{\mu_1}$ for model hypersurface of class $A_2$ in $\Bbb CH^m.$
Then $\mu_1, \; \mu_3$
have respective multiplicities $2l$ and $2k$ for some positive integers $k, \, l$ with $l = m - k - 1$ i.e. $n = 2l + 2k + 1.$ Moreover,
$$ \mu_1^* = \mu_1, \quad \mu_3^* = \mu_3, \quad \mu_1 \mu_3 = - c, \quad \k = \mu_1 - \frac{c}{\mu_1}= \mu_1 + \mu_3 \tag 59
$$
$$ f = L \mu_1 + K \mu_3, \quad f^2 = L^2\mu_1^2 + K^2 \mu_3^2 - 2c KL, \quad f_2 = L \mu_1^2 + K \mu_3^2 - 2c, \tag 60
$$
where $K:= 2k + 1$ and $L:= 2l + 1.$ Our goal is to examine when the equations ($E_1$)-($E_3$) are consistent and when constants $p$ and $q$
can be found to satisfy them (Once again, the condition ($E_4$) is satisfied by every class$-A_2$ hypersurface ). That comes down to the pair of
equations consisting of ($E_3$) and (50), having the same solution for $p$ for either value of $\mu \in \{\mu_1, \mu_3 \}.$
Consider the equation (50) in which $\mu = \mu_1,$
multiplied by $[2c + f (\mu_3 - \varkappa )] = (2c - f \mu_1)$ and the same equation with $\mu = \mu_3$ multiplied by $(2c - f \mu_3).$
Subtract the two multiplied equations to eliminate $p.$ We get
$$ f (f_2 + f^2) + 2 \k f (f + \k ) - c (n + 3) f - 4c \k = 0. \tag 61
$$
This is a necessary and sufficient condition for $p$ to have the same value from (50), regardless of the choice of $\mu .$ On the other hand subtracting the two equations obtained from (50) for $\mu = \mu_1, \, \mu_2,$ gives
$$ p f = f f_2 + c (3n + 13) f + 4c \k. \tag 62
$$
Similarly, from the two equations contained in ($E_3$) for $\mu = \mu_1, \; \mu_3$ by subtracting we get
$$ p = 2 f_2 + f^2 + 2 \k (f + \k) + 2c (n + 5), \tag 63
$$
and by eliminating $p$ from these two equations we get exactly the same condition  (61) as before. Moreover, assuming (61), we check that (62)
and (63) are consistent, so there is only one condition, namely (61), to be satisfied in order to make ($E_1$)- ($E_3$) consistent,
regardless of the choice of $\mu ,$ and enable us to solve for $p$ and $q$. Replacing the values from (59) and (60) into (61), using
$\varkappa = \mu_1 + \mu_3$ we get
$$ \align 0 = &\; L (L + 1)(L+2)\mu_1^3 + K (K+1)(K+2)\mu_3^3\\
& \; - c\mu_1 (3 L^2K + 3 L^2 + 6LK + 8L + 2K + 4 )\\
& \; - c \mu_3 (3 LK^2 + 3 K^2 + 6LK + 8K + 2L + 4),
\endalign
$$
or
$$ [(L+1)\mu_1^2 - c (K+1)]\,[ L(L+2)\mu_1^4 - 2c (LK + K + L + 2)\mu_1^2 + K(K+2)] = 0, \tag 64
$$
which has the following three solutions
$$ (a) \; \mu_1^2 = \frac{(K+1)c}{L + 1} \qquad\quad (b) \; \mu_1^2 = \frac{Kc}{L + 2} \qquad\quad
(c) \; \mu_1^2 = \frac{(K+2)c}{L}. \qquad\quad
$$
Clearly, when $c = - 1$ none of them is possible , so there are no 2-type hypersurfaces of $\Bbb CH^m(-4)$ among $A_2-$hypersurfaces.
When $c = 1$ the last two possibilities generate the same set of examples. From (63) we find
$$p = (L^2 + 4L + 2)\mu_1^2 + (K^2 + 4K + 2) \mu_3^2 - 2 LK$$
and we can also compute $q$ from ($E_1$) in terms of $\mu_1, \, \mu_3. $ Then using these we find the two eigenvalues of the Laplacian 
from the 2-type decomposition to be
$$\align \lambda_u &= (L+1)(L+2)\mu_1^2 + (K+1)(K+2)\mu_3^2 - (L + K + 2LK), \\
\lambda_v &= L\mu_1^2 + K\mu_3^2 + L + K, \qquad \; \mu_1 = \cot r, \; \mu_3 = - \tan r . \tag 65 \endalign $$
In the case $(a)$, we get $\lambda_u = 2(n + 3), \; \lambda_v = 2 (n+1) - \frac{(l-k)^2}{(l+1)(k+1)} 
= \frac{4(m+1)(LK+m)}{(L+1)(K+1)},$ 
$\; \lambda_u > \lambda_v, $ so the
hypersurface is of 2-type. Since $\mu_1^2 = \cot^2r = \frac{K+1}{L+1},$ from Takagi's list it follows that the hypersurface is an 
open portion of the tube of radius
$r = \cot^{-1}\sqrt{\frac{K+1}{L+1}} = \cot^{-1}\sqrt{\frac{k+1}{m-k}}$ about a canonically embedded
$\Bbb CP^k(4) \subset \Bbb CP^m(4),$ for any  $ k = 1,..., m-2;\,$ see also [8], [25]. For case $(b)$, (65) yields
$$ \lambda_u = \frac{4 (k + 1)(n+3)}{2k+1} = 4 (m+1)\frac{K+1}{K}, \quad  \lambda_v = \frac{4 (l + 1)(n+3)}{2l+3} 
= 4 (m+1) \frac{L+1}{L+2}, $$
\noindent
$\lambda_u > \lambda_v.$
Since $\mu_1^2 = \cot^2r = \frac{K}{L+2},$ we identify such hypersurface as an open portion of the tube of radius
$ r = \cot^{-1}\sqrt{\frac{K}{L+2}} = \cot^{-1}\sqrt{\frac{2k+1}{2(m-k)+1}}$ about a canonically embedded
$\Bbb CP^k(4) \subset \Bbb CP^m(4),$ for any  $ k = 1,..., m-2.$
\vskip 0.5truecm\noindent
$(ii)\;$ For an $A_2$-hypersurface we have from (65) 
$$ \lambda_u - \lambda_v = (L^2 + 2L + 2) \mu_1^2 + (K^2 + 2K + 2) \mu_3^2 - 2(L+K+LK). \tag{66}
$$
Note that from Tables 1 and 2 and the accompanying discussion, in addition to principal curvature $\k = 2 \cot_c(2r)$ 
an $A_2$-hypersurface has also two more principal curvatures $\mu_1 = \cot_c r$ and $\mu_3 = - c\tan_c r,$ 
with corresponding principal subspaces $V_1:= V_{\mu_1}$ and $V_3:= V_{\mu_3}, $ being $J$-invariant and $\; \Cal D = V_1 \oplus V_3.$
Then from (14) and (43) for a basis of principal directions $\{ e_i\}$ in $\Cal D$ we get
$$ \Delta \tx = - (L \mu_1 + K \mu_3) \xi - \s (\xi , \xi) - \sum_{e_i \in V_1} \s (e_i, e_i) - \sum_{e_j \in V_3} \s (e_j, e_j), 
$$ 
$$ \align \Delta^2 \tx = &- [L^2 \mu_1^3 + K^2 \mu_3^3 + (L^2 + 4mL - 4) \mu_1 + (K^2 + 4mK - 4) \mu_3]\, \xi\\
&+ [(L^2-2) \mu_1^2 + (K^2 - 2)\mu_3^2 - 2LK]\, \s (\xi , \xi)\\
&- 2 (L+1)(\mu_1^2 + 1) \sum_{e_i \in V_1} \s (e_i, e_i) - 2 (K+1)(\mu_3^2 + 1)\sum_{e_j \in V_3} \s (e_j, e_j).
\endalign
$$
Then $\tx_u $ and $\tx_v $ can be computed as in (55)-(56). Since the hypersurface of $\Bbb CQ^m$ is mass-symmetric via $\tx$
we have $\tx_0 = \tx - (\tx_u + \tx_v) = I/(m+1).$ Because $I$ and $\tx$ are normal to $\tx (\Bbb CQ^m),$ a necessary condition for 
mass-symmetry in $H^{(1)}(m+1)$ is that the $\xi$-component of $\tx_u + \tx_v$ be zero. The $\xi$-component of $\tx_u$ equals
$$ \frac{-4}{\lambda_u (\lambda_u - \lambda_v)} [(mL - 1)\mu_1 + (mK - 1)\mu_3]
$$
and the $\xi$-component of $\tx_v$ is
$$ \align \frac{1}{\lambda_v (\lambda_v - \lambda_u)} \{ L & (L^2 + 2L + 2) \mu_1^3  - [8mL + LK (3L + 2) + 2K - 4] \mu_1 \\
&  + K (K^2 + 2K + 2) \mu_3^3  - [8mK + LK (3K + 2) + 2L - 4]\, \mu_3  \}.
\endalign
$$
Observing the corresponding values of $\lambda_u , \, \lambda_v$ in each of the cases we see that the $\xi$-component of $\tx_u + \tx_v$
for hypersurfaces in $(b)$ is never zero, whereas for hypersurfaces of case $(a)$ this component is identically equal to zero. An additional 
computation verifies that for any hypersurface of case $(a)$ other components $\s (\xi , \xi), \, \sum_{e_i \in V} \s (e_i, e_i)$ 
on both sides of mass-symmetric 2-type decomposition are matched. 
\enddemo
\vskip 0.1truecm
The two families of tubes referred to in Lemma 3 have also another representation.
Let
$$M_{2k+1, 2l+ 1}(r):= S^{2k+1}(\cos r) \times S^{2l+1}(\sin r), \; 0 < r < \pi / 2, $$
be the family of generalized Clifford tori in an
odd-dimensional sphere $S^{n+2} \subset \Bbb C^{m+1}, \, n = 2m-1.$ By choosing the two spheres (with the indicated radii) in the above
product to lie in complex subspaces we get the fibration $S^1 \to M_{2k+1, 2l+1}(r) \to M^{\Bbb C}_{k, l}(r):= \pi (M_{2k+1, 2l+1}(r))$ 
compatible with the Hopf fibration $\pi : S^{n+2} \to \Bbb CP^m(4),$ which submerses $M_{2k+1, 2l+1}(r)$ onto $M^{\Bbb C}_{k, l}(r)$ [21].
Cecil and Ryan have shown [8] that $M^{\Bbb C}_{k, l}(r)$ is a tube of radius $r$ about totally geodesic $\Bbb CP^k(4)$
with principal curvatures $\cot r, \; - \tan r,\; 2 \cot (2r)$ of respective multiplicities $2l, 2k,$ and $1.$ Accordingly, the family
of hypersurfaces corresponding to the case $(a)$ is given as open portions of
\vskip 0.2truecm
$   \qquad\qquad   M^{\Bbb C}_{k, l}(r) = \pi \left(S^K\left( \sqrt{\frac{K+1}{n+3}}\right)\times
S^L\left( \sqrt{\frac{L+1}{n+3}}\right)  \right),\quad \cot^2r = \frac{K+1}{L+1},$
\newline\noindent
and the family of hypersurfaces corresponding to the case $(b)$ is
\vskip 0.2truecm
$ \qquad\qquad  M^{\Bbb C}_{k, l}(r) = \pi \left(S^K\left( \sqrt{\frac{K}{n+3}}\right)\times
S^L\left( \sqrt{\frac{L+2}{n+3}}\right)  \right), \quad \cot^2r = \frac{K}{L+2},$
\newline\noindent
where for both families $n+3 = 2 (m+1)$ and $K = 2k+1$ and $L = 2l+1$ are odd positive integers with $K + L = 2m.$ It is in
exactly this form that they appear in Udagawa's paper. The family of hypersurfaces corresponding to the case $(c)$ is the same family as in $(b),$
with the roles of $K$ and $L$ interchanged and the factors reversed. Hypersurfaces of case $(c)$ can be also described as tubes over 
$\Bbb CP^k(4)$ of radius $\rho = \cot^{-1}\sqrt{\frac{2k+3}{2(m-k)-1}},$ for $k = 1, 2, \dots , m-2, $ but are not listed as a separate
case since they constitute the same family as the one under case $(b)$. Namely, the tube about $\Bbb CP^k(4)$ of this radius $\rho$ is the same
as the tube over the other focal variety $\Bbb CP^l(4)$ of radius $\frac{\pi}{2} - \rho = \cot^{-1}\sqrt{\frac{2l+1}{2(m-l)+1}},$ which appears 
within family $(b).$
\vskip 0.4truecm\noindent
{\bf Remark.} Note that according to a result of Barbosa et al. [1], tubes over $\Bbb CP^k(4)$ of radius $r \,$ satisfying
$ \cot^{-1}\sqrt{\frac{2k+3}{2(m-k)-1}} \leq r \leq \cot^{-1}\sqrt{\frac{2k+1}{2(m-k)+1}} \;$ are stable with respect to normal
variations preserving the enclosed volume. Hence the 2-type tubes over $\Bbb CP^k(4)$ of radii $ \cot^{-1}\sqrt{\frac{2k+1}{2(m-k)+1}}$
and $ \cot^{-1}\sqrt{\frac{2k+3}{2(m-k)-1}}$
are distinguished by being maximal, respectively minimal, stable tubes over $\Bbb CP^k, $ for each  $ k = 1, 2, ..., m-2,$
i.e. the values of radii in cases $(b)$ and $(c)$ are precisely the endpoints of the stability interval for $r.$
\vskip 0.4truecm
\proclaim {Lemma 4} \,There are no 2-type hypersurfaces in $\Bbb CH^m(-4)$ among hypersurfaces of class $B.$ A class-$B$
hypersurface of $\Bbb CP^m(4)$ is of Chen 2-type if and only if it is an open portion of either the tube of radius
$r_1 = \cot^{-1}(\sqrt m + \sqrt{m+1})$ or the tube of radius $r_2 = \cot^{-1}\sqrt{\sqrt{2m^2-1} + \sqrt{2m^2-2}}, \; r_1 < r_2 ,$ 
about a complex
quadric $Q^{m-1} \subset \Bbb CP^m(4).$ In both instances, these tubes are also mass-symmetric in the hypersphere
$S^{N-1}_{I/(m+1)}\left(\sqrt{\frac{m}{2(m+1)}} \right)$ of $E^N = H(m+1)$ that contains them.
\endproclaim

\demo{Proof}
Let $\mu_2 = \cot (r + \frac{\pi}{4})$ and $\mu_4 = \cot (r + \frac{3\pi}{4}), \; \k = 2 \cot (2r)$ for the standard examples $B$ through
$E$ in
$\Bbb CP^m(4)$ and $\mu_2 = \coth r, \; \mu_4 = \tanh r , \; \k = 2 \tanh (2r)$ for an example of class $B$ in $\Bbb CH^m(-4).$
For all of these hypersurfaces we have
$$ \mu_2 \mu_4 = - c, \quad \mu_2 + \mu_4 = - \frac{4c}{\k }, \quad \mu_2^* = \mu_4, \; \; \mu_4^* = \mu_2. \tag 67
$$
Setting $\mu = \mu_2, \,\mu_4 $ in (50) produces two equations, from which by eliminating $p$ we get
$$ 2 [\k^2 + c (m-3) - \frac{16}{\k^2}]\, f + 2 (\k + \frac{4c}{\k}) f^2 - f (f^2 + f_2) - 4c\k = 0. \tag 68
$$
The same condition is obtained from ($E_3$) by setting $\mu = \mu_2,\, \mu_4 $ and eliminating $p$ and is also a necessary
condition for the values of $p$ obtained from (50) and ($E_3$) to be equal for any hypersurface of class $B.$
For class-$B$ hypersurface in either setting the common multiplicity of $\mu_2, \; \mu_4$ is $m-1$ and we have
$$ f = \k - \frac{4c (m-1)}{\k}, \qquad\quad f_2 = \k^2 + \frac{16 (m-1)}{\k^2} + 2c (m-1).\tag 69
$$
From (68) and (69) we get
$$ \k^6 - 4c (m-1)\k^4 - 8 (m^2 + 2m-1) \k^2 + 32c m (m^2 - 1) = 0.\tag 70
$$

Thus, for hypersurfaces of class $B,$ (70) represents a necessary and sufficient condition for $p$ to have the same value from (50) and ($E_3$),
regardless of the choice of $\mu = \mu_2, \, \mu_4,$ and also for $p$ and $q$ to be uniquely determined from the conditions
($E_1$)-($E_3$) in Lemma 1. One needs to check also condition ($E_4$) by computing the connection coefficients
of the hypersurface considered or by invoking the $\eta -$parallelisms of the shape operator for hypersurfaces of class $B, $ [20], [25].
We shall work instead with condition (44), which is a necessary and sufficient condition for type
$\leq 2,$ provided that the roots of the corresponding quadratic equation are real and distinct, and obtain expressions for
$\lp\tx$ and $\lp^2\tx$ that will enable us to find the explicit 2-type decomposition of $\tx$ for certain hypersurfaces of class $B.$
Let $M$ be such hypersurface in either $\Bbb CP^m(4)$ or $\Bbb CH^m(-4).$  With $\mu_2 , \mu_4$ as above, for the corresponding
eigenspaces $V:= V_{\mu_2}$ and $V_{\mu_4}$ we have $V_{\mu_4} = J V$ and $\Cal D = V \oplus SV.$
In the case of tube of radius $r = \frac 1 2 \ln (2 + \sqrt 3)$ in $\Bbb CH^m(-4)$ which has only two constant principal curvatures (since
$\varkappa = \mu_2),$ we consider $V_{\mu_2}$ to consists of eigenvectors of $\mu_2 $ belonging to $\Cal D $ only, thus not including $U.$
Note that
$\s (U, U) = \s (\xi , \xi)$ and $\s (Je_i, Je_i) = \s (e_i, e_i)$ by (12).

The fact that $J$ interchanges $V_{\mu_2}$ and $V_{\mu_4}$ for every hypersurface of class $B$ is crucial here and will enable us to find
suitable expressions for $\Delta \tx, \; \Delta^2\tx $ and later $\Delta^3 \tx .$ Let $\{e_i \} = \{e_j , Se_j\}$ be a $J-$basis of
the holomorphic distribution $\Cal D$ (where
$e_j \in V, \, j = 1, 2, ..., m-1$), which is the basis of principal directions of $A|_{\Cal D}$ with
$Ae_j = \mu_j e_j$ and $A(Se_j) = \mu^*_j Se_j,$ where $\mu_j, \mu_j^*$ satisfy relation (42), or equivalently
$$  2c + \k (\mu_j + \mu_j^*) = 2 \mu_j \mu_j^* = -2c. \tag 71
$$
Then we have
$$ \align \sum_{e_i\in\Cal D} \s (e_i, Ae_i) &= \sum_{e_i\in\Cal D} \mu_i \s (e_i, e_i) =
\sum_{e_j\in V}[\, \mu_j \s (e_j, e_j) + \mu_j^* \s (Se_j, Se_j)] \\
&= \sum_{e_j\in V} (\mu_j + \mu_j^*)\s (e_j, e_j) = - \frac{2c}{\k }\sum_{e_i\in\Cal D} \s (e_i, e_i), \quad \text{and}\tag 72
\endalign $$
$$\align \sum_{e_i\in\Cal D} \s (Ae_i, Ae_i) &= \sum_{e_i\in\Cal D} \mu_i ^2 \s (e_i, e_i) =
\sum_{e_j\in V} (\mu_j^2 + {\mu_j^*}^2)\s (e_j, e_j)\\
&= (\frac{8}{\k^2} + c) \sum_{e_i\in\Cal D} \s (e_i, e_i). \tag 73
\endalign
$$
Then using (67) and (69), formulas (14) and (18) become respectively

$$ \Delta \tx = \left[\frac{2c(n-1)}{\k} - \k \right] \xi - \s (\xi , \xi) - \sum_{e_i \in \Cal D} \s (e_i, e_i), \tag 74
$$
$$\align \Delta^2 \tx = & \left[ \frac{16c (n-1)^2}{\k^3} + \frac{8n(n-1)}{\k} - 2c (n+1)\k - \k^3  \right] \xi \\
&+ \left[ \frac{4(n-1)(n+3)}{\k^2} - 4c - \k^2 \right] \s (\xi , \xi ) - 2 (n+1)(\frac{4}{\k^2} + c )
\sum_{e_i \in \Cal D} \s (e_i, e_i). \tag{75}
\endalign
$$
For $X \in \Gamma(TM)$ we get the following using (3) and (11):
$$ \tilde\nabla_X \xi = - AX + \s (X, \xi ), \quad
\tilde\nabla_X (\s (\xi , \xi )) = - 2cX - 2c \la X, U\ra U - 2 \s (AX, \xi ), \tag 76
$$
$$ \tilde\nabla_X \left( \sum_{e_i \in \Cal D} \s (e_i, e_i )\right) = -2c(n+1)X
+ 4c \la X, U\ra U + 4 \s (AX, \xi ). \tag 77
$$
Therefore, differentiating $\lp\tx$ and $\lp^2\tx$ with respect to $X$ we substitute in (44) using (72)-(77) to get
$$\align &\left\{\frac{16c (n-1)^2}{\k^3} + \frac{2 (n-1)(4n-cp)}{\k} + [p - 2c (n+1)]\k - \k^3 \right\} [-AX + \s (X, \xi)]\\
& \qquad + \left[ \frac{8c (n^2 + 2n + 5)}{\k^2} + 4 (n^2 + 2n + 3) + 2c \k^2 + q - 2c (n+2)p\, \right] X \\
& \qquad + 2 \left[ p + \k^2 - 4cn - \frac{4(n^2 + 6n + 1)}{\k^2}\right] [c \la X, U\ra U + \s (AX, \xi)] = 0.\tag 78
\endalign
$$
Equate with zero the normal to $\Bbb CQ^m$ component of (78), which is a linear combination of $\s (X, \xi ) $ and $\s (AX , \xi),$
and take respectively $X \in V_{\mu_2}$ and $X \in V_{\mu_4}$ to get two equations, from which by subtracting and solving for $p$ we get
$$ p = \frac{4 (n^2 + 6n + 1)}{\k^2} + 4cn - \k^2.\tag 79
$$
Thus the last line of equation (78) drops out and the coefficient of $\s (X, \xi )$ on the top line would have to be zero.
With the value of $p$ from (79), equating that coefficient with 0 yields
$$ \k^6 - 2c (n-1)\k^4 - 2 (n^2 + 6n + 1)\k^2 + 4c (n-1)(n+1)(n+3) = 0. \tag 80
$$
Under this condition the $AX$- component is also zero and then the $X$- component in the middle line of (78) must be zero, which
gives the following value of $q:$
$$ q = 2c(n+3)\left[\frac{4(n^2 + 4n - 1)}{\k^2} + 2c (n-1) - \k^2 \right]. \tag 81
$$
Thus under the condition (80) it is possible to satisfy equation (78), that is the equation (44), for the values of $p$
and $q$ as in (79) and (81). This means that a class-$B$ hypersurface satisfying (80) is of 2-type. The equation (80)
is, not surprisingly, the compatibility condition (70), which we now see is also a sufficient condition for a hypersurface
of class $B$ to be of 2-type. Moreover, that condition is equivalent to the equation
$$ [\k^2 - 2c (n+1)]\, [\k^4 + 4c \k^2 - 2 (n-1)(n+3)] = 0, \tag 82
$$
which has three roots $\k^2 = 2c (n+1)$ and $\k^2 = - 2c \pm c \sqrt{2 (n^2 + 2n - 1)}.$ When $c = -1,$ none of them is possible
since $0 < \k^2 < 4.$ For $c = 1 $ (the case of a hypersurface of class $B$ in $\Bbb CP^m(4)$) we have the following two
possibilities:
$$ (a) \; \k^2 = 2 (n + 1) \qquad \text{and} \qquad (b) \; \k^2 = \sqrt{2 (n^2 + 2n - 1)} - 2.
$$
In case $(a)$ we find $p = \lambda_u + \lambda_v , \; q = \lambda_u \lambda_v $ from (79) and (81) and the two eigenvalues
$\lambda_u < \lambda_v$ to be
$$ p = \frac{4n(n+3)}{n+1}, \quad q = \frac{4 (n-1)(n+3)^2}{n+1} \tag 83
$$
$$ \lambda_u = \frac{2(n-1)(n+3)}{n+1} = 4 ( m -  1/ m ), \qquad \lambda_v = 2 (n+3) = 4 (m+1).\tag 84
$$
The corresponding hypersurface is a tube of radius $r$ about the complex quadric $Q^{m-1},$ where $\cot r - \tan r = \k, $
i.e. $\cot r = \frac{\k + \sqrt{\k^2 + 4}}{2} = \sqrt m + \sqrt{m+1}.$

In case $(b)$ we find
$$ p = \frac{1}{(n-1)(n+3)}\, [\,2 (2n^3 + 7n^2 + 8n - 1) + (n^2 + 10n + 5) \sqrt{2 (n^2 + 2n - 1)}\,], \tag 85
$$
$$ q = \frac{2}{n-1}\, [\, 2 (n^3 + 4n^2 + 5n - 2 ) + (n^2 + 6n + 1)\sqrt{2 (n^2 + 2n - 1)}\, ], \tag 86
$$
$$ \lambda_u = \frac{2}{n+3}\, [2 \sqrt{2 (n^2 + 2n - 1)} + (n+1)^2], \quad
\lambda_v = \frac{n+3}{n-1}\, [\sqrt{2 (n^2 + 2n - 1)} + 2n ].\tag 87
$$
The corresponding hypersurface is a tube about $Q^{m-1}$ of radius $r,$ with $\cot r = \sqrt{\sqrt{2m^2-1} + \sqrt{2m^2-2}}.$
These two tubes are therefore of 2-type in $H(m+1).$ Moreover, they are also mass-symmetric in the hypersphere containing
$\Phi (\Bbb CP^m ),$ which means that the center of mass is $\tx_0 = I/(m+1) = 2I/(n+3).$ Indeed by a Lemma of [15]
we have the expression
$$ I = (m+1) \tx + \frac{c}{2}\, \s (\xi , \xi ) + \frac{c}{4} \sum_{e_i \in \Cal D} \s (e_i, e_i).\tag 88
$$
Then it is a straightforward verification using (74), (75), (79), (81) and (88) that any class-$B$ hypersurface satisfying
condition (80) is of 2-type since it satisfies the equation
$$ \Delta^2 \tx - p\, \Delta\tx + q \left( \tx - \frac{I}{m+1}\right) = 0, \tag 89
$$
and it is, obviously, not of 1-type. Specifically, for the two tubes about $Q^{m-1}$ discussed above, (89) holds for the indicated
values of $p$ and $q$ from (83), respectively (85)-(86). The corresponding vector-eigenfunctions $\tx_u $ and $\tx_v $ of
$\lambda_u $ and $\lambda_v $ in 2-type decomposition of $\tx $ can be found from
$$ \tx_u = \frac{1}{\lambda_u - \lambda_v }\left[\Delta \tx - \lambda_v \left(\tx - \frac{I}{m+1}\right)\right], \;\;
\tx_v = \frac{1}{\lambda_v - \lambda_u }\left[ \Delta \tx - \lambda_u \left(\tx - \frac{I}{m+1}\right)\right]. \tag 90
$$

For example, for the tube of radius $r_1 = \cot^{-1}(\sqrt m + \sqrt{m+1}) $ we get
$$ \tx_u = \frac{\sqrt{2(n+1)}}{2(n+3)}\, \xi - \frac{n+1}{4 (n+3)}\, \s (\xi , \xi ), \quad \text{and }
$$
$$ \tx_v = - \frac{\sqrt{2(n+1)}}{2(n+3)}\, \xi + \frac{n - 3}{4 (n+3)}\, \s (\xi , \xi )
- \frac{1}{2(n+3)}\sum_{e_i \in \Cal D} \s (e_i, e_i).
$$
It can be also directly verified, using (16), (26) and (17), that $\tx_u, \, \tx_v$ are indeed eigenfunctions of $\Delta$ 
for the indicated eigenvalues.
Incidental to this finding, we obtain two simple eigenvalue estimates for the first two non-zero eigenvalues
$\lambda_1, \, \lambda_2$ for the hypersurface for which $(a)$ holds: $\, \lambda_1 \leq 4 (m - 1/m)$ and $\lambda_2 \leq 4 (m+1).$ 
\enddemo
\vskip 0.1truecm

\proclaim {Lemma 5} \, There are no 2-type hypersurfaces in $\Bbb CP^m (4)$ among any of the standard examples of class
$ C,\, D, \,$ or $E.$
\endproclaim
\demo{Proof} This was shown in [36]. For the sake of completeness we include a different proof here using our approach.
In addition to principal curvatures $\mu_2, \mu_4$ and formulas (67)-(68), we have also principal curvatures $\mu_1, \mu_3,$
for which the relations (59) hold ($c = 1$ throughout). If we substitute $\mu = \mu_1,\, \mu = \mu_3$ in ($E_3$) and subtract the
two resulting equations we get
$$ p = 2 f_2 + f^2 + 2\k f + 2 \k^2 + 2 (n+5). \tag 91
$$
The same manipulation with $\mu = \mu_2, \, \mu_4$ yields
$$ p = 2 f_2 + f^2 - 2 (\k + \frac{4}{\k})f + 2 (n+5) + \frac{32}{\k^2} - 2 \k^2. \tag 92
$$
On the other hand, substituting $\mu = \mu_1, \mu_3$ into (50) and subtracting we get
$$ fp = f f_2 + (3n + 13)f + 4\k, \tag 93
$$
and the same procedure using  $\mu = \mu_2, \mu_4$ leads to
$$ fp = f f_2 + (3n+5)f - 4 \k . \tag 94
$$
Combining (93) and (94) we get $f = - \k $ and subtracting (91) and (92) leads to $(\k + 2/\k)f + \k^2 - 8/\k^2 = 0,$
which is incompatible with $f = - \k.$ 
\enddemo
\vskip 0.1truecm

Now we can formulate our main classification results for 2-type Hopf hypersurfaces of $\Bbb CQ^m.$
In the complex projective space we have

\proclaim{Theorem 1} Let  $M^{2m-1}$ be a Hopf hypersurface of $\Bbb CP^m(4),
 \, ( m \geq 2)$. Then $M^{2m-1}$ is of 2-type in $H(m+1)$ via $\Phi$ if and only if it is an
 open portion of one of the following
 \roster
 \item "{\it (i)}" A geodesic hypersphere of any radius $r \in (0, \frac{\pi}{2}),$
 except for  $r = \cot^{-1}\sqrt{\frac{1}{2m+1}};$
 \item "{\it (ii)}" The tube of radius $r = \cot^{-1}\sqrt{\frac {k+1}{m-k}}$ about a canonically embedded totally geodesic
$\Bbb CP^k(4) \subset \Bbb CP^m(4),$ for any $k= 1,..., m-2;$
\item "{\it (iii)}"  The tube of radius $r = \cot^{-1}\sqrt{\frac{2k+1}{2(m-k)+ 1}}$  about a canonically embedded
$\Bbb CP^k(4) \subset \Bbb CP^m(4),$ for any $ \, k = 1, ..., m-2.$
\item "{\it (iv)}" The tube of radius $r = \cot^{-1} (\sqrt m + \sqrt{m+1})$ about a complex quadric $Q^{m-1}\subset \Bbb CP^m(4).$
\item "{\it (v)}" The tube of radius $r = \cot^{-1}\sqrt{\sqrt{2m^2-1} +  \sqrt{2m^2-2}}$
about a complex quadric $Q^{m-1}\subset \Bbb CP^m(4).$

\endroster
 \endproclaim
\demo{Proof} As shown  before, 2-type Hopf hypersurface must have constant principal curvatures and therefore must be one from
the Takagi's list in $\Bbb CP^m(4).$ The rest follows from Lemmas 1-5.
\enddemo
\vskip 0.1truecm
As commented before, the same classification holds when $M$ is assumed to have constant mean curvature (CMC) instead of being Hopf.
In that regard Theorems 1 and 2 in [36] are deficient and incomplete since Udagawa's list contains examples $(i)-(iii)$ only. The list of 
items $(i)-(v)$ is the correct and complete classification of CMC hypersurfaces of 2-type in $\Bbb CP^m(4).$

\par
In the same manner, since being a Hopf hypersurface and having constant mean curvature imply each other for hypersurfaces of 2-type,
Lemmas 1-5 yield

\proclaim{Theorem 2}\, Let  $M^{2m-1}$ be a real hypersurface of $\Bbb CH^m(-4),
 \, ( m \geq 2)$ for which we assume that it is a Hopf hypersurface or has constant mean curvature. Then $M^{2m-1}$ is of 2-type
in $H^1(m+1)$ via $\Phi$ if and only if it is (an open portion of) either a geodesic hypersphere of arbitrary radius $r > 0$ or a tube of
arbitrary radius $r > 0$ about a canonically embedded totally geodesic complex hyperbolic hyperplane $\Bbb CH^{m-1}(-4).$
\endproclaim

Regarding mass-symmetric hypersurfaces, from the analysis above we have

\proclaim {Corollary 1} A complete Hopf (or CMC) hypersurface of $\Bbb CP^m(4)$ is of 2-type and mass-symmetric in the hypersphere of $H(m+1)$ 
containing $\Phi (\Bbb CP^m)$ if and only if it is one of the hypersurfaces (tubes) in $(ii), (iv) $ and $(v)$ or the geodesic hypersphere 
of radius $\cot^{-1}(1/\sqrt m).$ There exists no 2-type mass-symmetric (in particular, no null 2-type) hypersurface of $\Bbb CH^m(-4).$
\endproclaim

This rectifies the claim made in Theorem 2 of [36].

\head 6. CMC Hopf Hypersurfaces of 3-Type
\endhead
\vskip 0.4truecm

It is not difficult to see that the hypersurfaces of class $A_2$ are, generally speaking, of 3-type (except for those two families of tubes
in $\Bbb CP^m$ given in Theorem 1 $(ii), (iii),$ which are or 2-type).
Consider $p \in M \subset \Bbb CQ^m(4c)$ where $p = [\zeta ]$ is represented by a column vector
$$\zeta \in \pi^{-1}(p) \subset N^{2m+1} \subset \Bbb C^{m+1}_{(1)} = \Bbb C^{k+1}_{(1)}\oplus \Bbb C^{l+1}.$$
Let $z = (z_i) = (\zeta_0, ..., \zeta_{k})^T$ and $w = (w_\alpha) = (\zeta_{k+1}, ..., \zeta_{m})^T$  and consider in $\Bbb C^{k+1}_{(1)}$ 
the quadric $N^{2k+1}(r_1)$ (the sphere or anti - de Sitter space of radius $r_1$) and in $\Bbb C^{l+1}$ the sphere $S^{2l+1}(r_2)$ so 
that $r_1^2 + cr_2^2 = 1.$ In the projective case we have $c = 1$ and we set $r_1 = \cos r, \; r_2 = \sin r,$ whereas in the hyperbolic 
case $c = -1, \; r_1 = \cosh r, \; r_2 = \sinh r.$ The corresponding class$-A_2$ hypersurfaces which are the tubes of radius $r$ about 
totally geodesic $\Bbb CQ^k(4c)$ are obtained as the Hopf projections, defining the submersion:
$\pi \left( S^{2k+1}(\cos r) \times S^{2l+1}(\sin r)\right)$ in $\Bbb CP^m(4)$ and $\pi\left( H^{2k+1}_1(\cosh r) \times S^{2l+1}(\sinh r) \right) $ in
$\Bbb CH^m(-4), \; k + l = m-1.$ According to (2), the coordinate representation of $\tx (p)$ in $H^{(1)}(m+1)$ has the matrix block form
$$ \tx = \pmatrix a_{ij} & b_{i\beta}\\ c_{\alpha j} & d_{\alpha\beta}
\endpmatrix , \quad 0 \leq i,j \leq k, \;\; k+1 \leq \alpha , \beta \leq m,
$$
where, for example, $d_{\alpha\beta} = c\, w \bar w^T, \; b_{i\beta} = c\, z \bar w^T, $\ and $a_{ij} = (\pm z_i \bar z_j)$ is formed by
the signed products, plus in the first column minus otherwise in $\Bbb CH^m$-case, all plus in $\Bbb CP^m$-case.
Then using the fact that $\pi$ is a (pseudo) Riemannian submersion with totally geodesic fibers [2], one can compute the iterated
Laplacians of $ \pi (N^{2k+1} (r_1) \times S^{2l+1} (r_2))$ as follows, see [23], [36], [18]:

$$ \lp \tx = \pmatrix \frac{2c(K+1)}{r_1^2} a_{ij} - 4c I_{k+1} & \left(\frac{cK}{r_1^2} + \frac{L}{r_2^2}\right) b_{i\beta }\\
\left(\frac{cK}{r_1^2} + \frac{L}{r_2^2}\right) c_{\alpha j}  &  \frac{2(L+1)}{r_2^2} d_{\alpha \beta} - 4c I_{l+1}
\endpmatrix ,
$$
and in general for an integer $s \geq 1$
$$ \lp^s \tx = \pmatrix \frac{2^sc^s(K+1)^s}{r_1^{2s}} a_{ij} - \frac{2^{s+1}c^s(K+1)^{s-1}}{r_1^{2(s-1)}} I_{k+1}
& \left(\frac{cK}{r_1^2} + \frac{L}{r_2^2}\right)^s b_{i\beta} \\
\left(\frac{cK}{r_1^2} + \frac{L}{r_2^2}\right)^s c_{\alpha j}
&  \frac{2^s(L+1)^s}{r_2^{2s}} d_{\alpha \beta} - \frac{2^{s+1}c(L+1)^{s-1}}{r_2^{2(s-1)}} I_{l+1}
\endpmatrix .
$$

Then one checks that the following equation is satisfied
$$ \lp^3\tx + p \lp^2 \tx + q \lp \tx + r (\tx - \tx_0) = 0 \tag 95
$$
for
$$p = - \left[\frac{c(3K+2)}{r_1^2} + \frac{3L+2}{r_2^2} \right], \quad
r = - \frac {4c(K+1)(L+1)}{r_1^2 \, r_2^2} \left( \frac{cK}{r_1^2} + \frac{L}{r_2^2}\right),$$
$$ q = 2 \left [  \frac{K (K+1)}{r_1^4} + \frac{L (L+1)}{r_2^4} + \frac{c\, (4KL + 3K + 3L + 2)}{r_1^2\, r_2^2}\right ],
$$
and
$$ \tx_0 = \pmatrix \frac{2r_1^2}{K+1} I_{k+1} & O\\O & \frac{2cr_2^2}{L+1} I_{l+1} \endpmatrix , \quad k+l = m-1. \tag 96
$$
This means that any $A_2-$hypersurface is of 3-type if the polynomial
$\lambda^3 + p\, \lambda ^2 + q\, \lambda + r $ has simple real roots (and the hypersurface is not already of lower type).
Those roots are found to be
$$ \lambda_u = \frac{cK}{r_1^2} + \frac{L}{r_2^2}, \qquad \lambda_v = \frac{2c(K+1)}{r_1^2}, \qquad \lambda_w = \frac{2(L+1)}{r_2^2}.
$$

When $c = 1,$ the equality of any two among these three roots leads to 2-type examples $(ii)$ and $(iii)$ in Theorem 1.
If we look for mass-symmetric examples of class $A_2$ then $\tx_0 = I/(m+1),$ which gives $\cot^2r = \frac{K+1}{L+1},$ 
thus again leading to the example $(ii),$ which is of 2-type. So there are  no mass-symmetric 3-type examples among $A_2-$hypersurfaces 
in $\Bbb CP^m.$ On the other hand, when $c = -1 $ no equality between the roots $\lambda_u, \, \lambda_v, \,\lambda_w \, $ is possible 
and we know from Lemma 3 that no example of class $A_2$ in $\Bbb CH^m(-4)$ is of 2-type, they are all, therefore, of 3-type. Since the 
constant part $\tx_0$ in 3-type decomposition has the form given in (96) and cannot clearly equal $I/(m+1),$ the only way such hypersurface 
can be mass-symmetric, according to our definition,  is that the hypersurface is of null 3-type, i.e. the eigenvalue $\lambda_u = 0,$ in 
which case $\tx_0$ can be changed 
to equal $I/(m+1).$ This gives the condition $\coth^2r = K/L, $ i.e. the radius of the tube about $\Bbb CH^k(-4)$ is 
$r = \coth^{-1}\sqrt{\frac{2k+1}{2l+1}}, \; 1 \leq l < k \leq m-2, \; k + l = m - 1.$ In that case we get a mass-symmetric null 3-type 
hypersurface in $\Bbb CH^m(-4):$
$$ \pi (H^K_1(\cosh r) \times S^L(\sinh r)) = \pi \left( H_1^{2k+1}\left( \sqrt{\frac{2k+1}{2(k-l)}}\right) \times 
S^{2l+1}\left( \sqrt{\frac{2l+1}{2(k-l)}}\right)\right).
$$
Additional examples of mass-symmetric 3-type hypersurfaces have to be searched for among classes $B,\, C,\, D,\, $ and $E.$
We derive next certain necessary conditions for hypersurface with $\text{tr}\, A =$ const to be mass-symmetric and of 3-type.

\vskip 0.3truecm
Let $M^n$ be a CMC Hopf hypersurface of $\cxf (4c), \, (n = 2m-1)$ which is of 3-type via $\tx$ and mass-symmetric in the hyperquadric
centered at $I/(m+1)$ containing $\Phi (\cxf )$ and defined by $\langle P -\frac{I}{m+1}, P - \frac{I}{m+1} \rangle = \frac{cm}{2(m+1)}.$ Then
$$ \lp^3\tx + p \lp^2 \tx + q \lp \tx + r (\tx - I/(m+1)) = 0, \tag 97
$$
where $p, q, r $ are the (signed) elementary symmetric functions of the eigenvalues $\lambda_u, \lambda_v, \lambda_w $ associated with 
a 3-type decomposition of $\tx$. We will consider various components of this equation  to derive a set of necessary conditions for a
Hopf hypersurface with constant $\text{tr}\, A $ to be mass-symmetric and of 3-type. Those will include the conditions
$\text{tr}\, A^k = \text{const}, \, 1 \leq k \leq 4.$ Recall that the normal space $T^\perp_P\Bbb CQ^m $ in $H^{(1)}(m + 1)$ is spanned by
the position vector $P$ and vectors of the form $\s (Z, W), \, Z, W \in T_P\Bbb CQ^m $ [15]. Using (12) and (23)
we get from (14), (18), and (34) respectively
$$ \la \lp \tx , \tx \ra = n, \qquad \la \lp^2 \tx , \tx \ra = f^2 + 2c (n^2 + 2n - 1), \tag{98}
$$
$$ \align \la \lp^3 \tx , \tx \ra = \; f^2&[f_2 + 5c(n-1)] - 8c \k^2 + 4n (n+2)^2 - 20 \\
& + 16c\k f - 16 c f\, \text{tr}\, (SAS) - 4c\, \text{tr}\, (SA^2S) + 8c\, \text{tr}\, (SA)^2. \tag 99
\endalign
$$
Further, choosing a $J-$basis $\{ e_i, Se_i \}$ of $\Cal D$ and using (42) in the form
\newline
$\mu_i\mu_i^* = c + \frac{\k}{2}(\mu_i + \mu_i^*),$
we compute
$$ \text{tr}\, (SAS) = \k - f, \quad \text{tr}\, (SA^2S) = \k^2 - f_2, \quad \text{tr}\, (SA)^2 = \k^2 - \k f - (n-1)c.
$$
Substituting in (99) we obtain $\la \lp^3\tx , \tx\ra $ as a sum of several terms, one of which is $(f^2 + 4c)f_2$ and the others are
constants depending only on $\k , f, c, n.$ Therefore taking the metric product of (97)
with $\tx$ and using the above information, we see that if $f^2 + 4c \neq 0 $ (this condition is always satisfied in the projective case) 
it follows
that $f_2 = \text{tr}\, A^2$ is constant. Thus we will subsequently assume that $f^2 \neq 4 $ in the hyperbolic case to ensure the
constancy of $f_2.$

Next, we look at the $\xi-$component of (97). From (5) we compute
$$ \sum_i \la J [(\nabla_{e_i}A)^2(Se_i) - (\nabla_{e_i}A)(SAe_i)], \; \xi \ra = \k^2 f - \k f_2,
$$
so that the $\xi-$component of $\lp^3\tx $ equals
$$ \align \la \lp^3\tx , \xi \ra = 8&c (\k^2 f - \k f_2 ) + 8c\k^3 + 8 (2cf_2 + n + 4) \k\\
& - 8cf_3 - f [f^2_2 + 4c (n+4)f_2 + 4cf^2 + 7n^2 + 30 n + 19 ].\tag{100}
\endalign
$$
We also have that $\la \lp\tx , \xi \ra = - f $ and $\la \lp^2 \tx, \xi \ra = 4c\k - f [f_2 + c (3n+5)]$ are constant. Thus taking the metric
product of (97) with $\xi$ we get $f_3 = $ const . In finding the $\s (\xi , \xi)-$component
of (97) note that
$$ \la \s (\xi , \xi ), \s (\xi , \xi )\ra = 4c , \quad \sum_i \la \s (e_i, SAS e_i), \s (\xi , \xi )\ra = -2c (f - \k),
$$
$$ \sum_i \la \s (e_i, SA^2Se_i), \s (\xi , \xi ) \ra = - 2c (f_2 - \k^2),
$$
$$
\sum_i \la \s (e_i, (SA)^2e_i), \s (\xi , \xi ) \ra = - 2 (n-1) - 2c\k (f - \k ),
$$
and
$$ \align \sum_{i, j} &\la \s ((\nabla_{e_j}A)e_i, (\nabla_{e_j}A)e_i ), \s (\xi , \xi)\ra\\
&= 2 c \Vert \nabla A\Vert^2 + 2c \sum_j \la (\nabla_{e_j}A)U, (\nabla_{e_j}A)U \ra\\
&= 2 c \Vert \nabla A\Vert^2 - 2c \k^2 \, \text{tr}\, (AS^2A) + 4c\k \, \text{tr}\, (ASASA) - 2c\, \text{tr}\, (ASA)^2 \\
&= 2c [f_2 - c (n+3)]\, f_2 - 2c f f_3 + 14(n-1)c\\
&\qquad + 2f^2 + 6 \k f + c \k^2 f_2 - (n-1)\k^2 - c \k^3 f ,
\endalign
$$
by way of (23) and (71). Now, when the inner product of (97) with $\s (\xi , \xi )$ is taken we get a sum of several terms equal to zero. The
only term in this sum containing $f_4$ is $8c f_4 $ and the other terms depend only on $\k , f, f_2 , f_3, p, q, r, c, n, $ and are thus constant.
It follows, therefore, that $f_4 = $\, const. We note that $\s (U, \xi ) = 0 $ by (12) and
considering $\s (X, \xi )-$component of (97)  for  $X \in \Gamma (\Cal D )$ we compute using the Codazzi equation that $\text{tr}\, 
((\nabla_XA)\circ [A, S]) = 0.$ For the part of (97) normal to $\cxf (4c)$ it remains to consider
$\s (X, Y)-$component for $X \in \Gamma (\Cal D )$ and $Y \in \Gamma (TM). $
We compute
$$ \la \s (\xi , \xi ), \s (X, Y )\ra = 2c \la X, Y\ra + 2c \la X, U\ra \la Y, U \ra,
$$
$$ \sum_i \la \s (A^ke_i , A^l e_i), \s (X, Y)\ra = 2c\, [\text{tr}\, (A^{k+l})\la X, Y\ra + \la A^{k+l}X, Y\ra - \la SA^{k+l}SX, Y \ra ],
$$
with $k, l$ integers $\geq 0, $ and $A^0 = I.$ Further,
$$ \sum_i \la \s (e_i, SASe_i), \s (X, Y)\ra = 2c (\k - f)\la X, Y \ra - 2c \la AX, Y\ra + 2c \la SASX, Y\ra ,
$$
$$\sum_i \la \s (e_i, SA^2Se_i), \s (X, Y)\ra = 2c (\k^2 - f_2)\la X, Y \ra - 2c \la A^2X, Y\ra + 2c \la SA^2SX, Y\ra ,
$$
$$\align \sum_i \la \s (Ae_i, &SASe_i), \s (X, Y)\ra = \sum_i \la \s (e_i, SASAe_i), \s (X, Y)\ra\\
&= 2c\, [\k^2 - \k f - (n-1)c]\la X, Y \ra + 2c\, \la [(SA)^2 + (AS)^2]X,\, Y\ra ,
\endalign
$$
$$ \align \sum_{i, j} &\la \s ((\nabla_{e_j}A)e_i, (\nabla_{e_j}A)e_i ), \s (X, Y)\ra\\
&= 2 c\, \Vert \nabla A\Vert^2 \la X, Y\ra + 2c\, \la BX, Y\ra + 2c\, \la B(SX), SY \ra ,\endalign
$$
where $B:= \sum_j (\nabla_{e_j}A)^2 $ is a well-defined endomorphism of $TM,$ independent of the choice of the
basis $\{e_i \}.$ Next, we compute
$$\la \lp \tx , \s (X, Y)\ra = - 2c (n + 2)\la X, Y\ra , $$
$$ \align \la \lp^2 \tx , \s (X, Y)\ra = \; &-4c \la ASX, ASY\ra - 4cf \la ASX, SY\ra - 4c \la AX, AY \ra \\
&- 4cf \la AX, Y\ra - [4 (n + 1)(n + 3) + 2c f^2 ] \la X, Y\ra,
\endalign $$
$$ \align\la \lp^3 \tx , \s (X, Y)\ra = \; &8c \la BX, Y\ra + 8c \la BSX, SY \ra + 32 \la [(SA)^2 + (AS)^2]X, Y \ra \\
&+ 8f (cf_2 + n + 7) \la SASX, Y \ra  + 16 (1 + c f_2) \la SA^2SX, Y\ra \\
&+ 8c \la SA^4SX, Y\ra - 8c \la A^4X, Y \ra - 16 (1 + cf_2) \la A^2X, Y\ra\\
&- 8f (cf_2 + n + 7)\la AX, Y\ra + [8 \k^2 + 16 \k f - 8 f_2 - 2c f^2 f_2 \\
&\qquad\quad - 2 (5n + 19)f^2 - 8c (n^3 + 6n^2 + 10n + 7)] \la X, Y\ra . 
\endalign
$$
Thus taking the inner product of (97) with $\s (X, Y)$ and dropping $Y$ we get
$$\align BX - SBSX = \; &A^4X - SA^4SX + a (A^2X - SA^2SX)\\
&+ b (AX - SASX) - 4c\, [(SA)^2 + (AS)^2] X + d X, \tag{101}
\endalign
$$
where $ X \in \Gamma (\Cal D), \; a = (p/2) + 2 (c + f_2), \; $ $b = (pf/2) + cf (cf_2 + n + 7), $ and
$$\align d = n^3 &+ 6n^2 + 10n + 7 + \frac{c}{4}(5n + 19) f^2 + \frac{1}{4}f^2 f_2 + c f_2 - 2 c \k f - c \k^2\\
& + \frac{p}{4} [2c (n+1)(n+3) + f^2] + \frac{q}{4} (n+2) + \frac{cr}{8}.
\endalign $$
There remains the part of (97) tangent to $M$ to be considered. Relation (97) has no $U-$component and for $X \in \Cal D$ we have by the Codazzi equation
$$ \sum_i \la J[(\nabla_{e_i}A^2)(Se_i) - (\nabla_{e_i}A)(SAe_i)],\, X \ra = \text{tr}\, ((\nabla_{SX}A)\circ [S, A]).
$$

The right-hand side of this is also the result of the metric product of the left-hand side of (97) with $X$ and therefore must be equal
to zero, which is the same piece of information contained in the $\s (X, \xi )-$component. Hence we have the following
\proclaim{Lemma 6} Let $M^n \subset \cxf (4c)$ be a Hopf hypersurface (n = 2m-1) with constant mean curvature (in the hyperbolic case we 
assume, additionally, that $\, (\text{tr}\, A)^2 \neq 4$). If $M^n$ is mass-symmetric and of 3-type in $H^{(1)}(m + 1)$ then we have
\roster
\item "{($i$)}" $\text{tr}\, A^k = \text{const}$, for $k = 1, 2, 3, 4;$
\item "{($ii$)}" $\text{tr}\, ((\nabla_X A) \circ [A, S]) = 0,$ for every $X \in \Gamma (\Cal D );$
\item "{($iii$)}" $$\align BX - SBSX = \; &A^4X - SA^4SX + a (A^2X - SA^2SX)\\
&+ b (AX - SASX) - 4c\, [(SA)^2 + (AS)^2] X + d X, \endalign
$$
where $B := \sum_j (\nabla_{e_j}A)^2, \; X \in \Gamma (\Cal D ), \; c = \pm 1, $ and $a, b, d$ are constants.
\endroster
\endproclaim
Essentially, the conditions $(i)- (iii)$ are also sufficient conditions for such $M$ to be mass-symmetric and of type $\leq 3$ since we
obtained these conditions by considering all of the components of equation (97), provided that the constants $a, b,c, d$ and $\text{tr}\, A^k$
are such to enable $p, q, r $ to be real and the polynomial equation $t^3 + p t^2 + q t + r = 0 $ to have simple roots [12].
Note that the condition $(ii)$ is automatically satisfied when $M$ is a Hopf hypersurface whose induced Hopf foliation is a Riemannian 
foliation [3, p. 64] since then $U$ is a Killing vector field and $AS = SA.$ See also [27], characterizing class$-A$ hypersurfaces by the condition $AS = SA.$

\proclaim{Corollary 2} Let $M$  be a CMC Hopf hypersurface of $\cxf $ satisfying $(\text{tr}\, A)^2 \neq - 4c $ and having at most four distinct
principal curvatures at each point. If $M$ is mass-symmetric and of 3-type via $\tx$ then $M$ has constant principal curvatures.
\endproclaim

We show next that every hypersurface of class $B$ is mass-symmetric in the corresponding hyperquadric and,
apart from those two tubes in Lemma 4, of 3-type.
Using the information from Lemma 4, we get respectively from (16), (26), (17) and (67) the following
$$ \Delta \xi = [f_2 + c (n-1)] \xi + (2\k - f) \s (\xi , \xi) - \frac{4c}{\k} \sum_{e_i\in\Cal D} \s (e_i, e_i), \tag 102
$$
$$ \Delta (\s (\xi , \xi)) = 4c \k \xi + 2 \left[ \frac{8(n-1)}{\k^2} + c (n+1)\right]\, \s (\xi , \xi) - \frac{16}{\k^2} 
\sum_{e_i\in\Cal D} \s (e_i, e_i), \tag 103
$$
$$ \align \sum_{e_i \in \Cal D} \Delta (\s (e_i, e_i)) = &\; 2 \left[ \frac{16}{\k^2} + c (n+3)\right] \sum_{e_i\in\Cal D} \s (e_i, e_i)
- 2 \left[\frac{16(n-1)}{\k^2} + c (n-1) \right] \s (\xi , \xi )\\
& + 2c \left[ (n-1)\k - \frac{2c(n-1)(n+3)}{\k} \right] \xi . \tag 104
\endalign $$

The third iterated Laplacian $\Delta^3\tx $ for a hypersurface of class $B$ can be computed from (34)-(35) but that would require
finding the connection coefficients of the hypersurface. It seems easier to find $\Delta^3 \tx$ directly by applying the Laplacian
to (75) and using (102)-(104). We get
$$ \align \Delta^3 \tx = &\; [ \frac{128c (n-1)^3}{\k^5} + \frac{128 (n-1)(n^2 + 1)}{\k^3} + \frac{8c (n-1)(3n^2 + 2n + 3)}{\k} \\
&\qquad\qquad\qquad\qquad\qquad - 16n \k - 4c (n+1)\k^3 - \k^5 ] \, \xi \\
& + [ \frac{32(n-1)(n+3)(3n+1)}{\k^4} + \frac{8c (n-1)(3n^2 + 14n + 3)}{\k^2} \\
&\qquad\qquad\qquad\qquad + 8 (n^2 - 4n + 1) - 2c (3n+1)\k^2 - \k^4 ] \, \s (\xi , \xi)\\
&- [ \frac{128 (n+1)^2}{\k^4} + \frac{48c (n+1)^2}{\k^2} + 4 (n-1)(n+3) - 4c \k^2 ] \, \sum_{e_i \in \Cal D} \s (e_i, e_i).
\tag 105 \endalign
$$
Then using (74), (75), and (105) we have
\proclaim {Lemma 7} Every class-$B$ Hopf hypersurface with constant principal curvatures in  $\Bbb CH^m(-4)$ is
mass-symmetric and of 3-type via $\tx .$ The same is true for hypersurfaces of class $B$ in $\Bbb CP^m(4), $ with
the exception of those two tubes about $Q^{m-1}$ referred to in Lemma 4, which are mass-symmetric and of 2-type.
\endproclaim
\demo{Proof} These hypersurfaces are not of 1-type and the only 2-type examples are given in Lemma 4. We show that they satisfy 
the 3-type equation (97) and that the polynomial $P(t) = t^3 + p t^2 + q t + r $ has three distinct real roots except for one value 
of $\k ,$ so that for other values of $\k$ the result of [12] proves it then to be
mass-symmetric and of type $\leq 3.$ Moreover, they will be exactly of 3-type if $P(t)$ is the minimal polynomial of the immersion 
$\tx - \tx_0, $
i.e. the hypersurface does not satisfy  a lower-degree polynomial in the Laplacian. Indeed, using the Gauss elimination with
$$ p = - \frac{1}{\k^2}(\k^2 + 4c)[\k^2 + 2c (3n+1)], $$
$$ q = \frac{4}{\k^4}(\k^2 + 4c) [c (n+1) \k^4 + (3n^2 + 6n - 1) \k^2 + 8c (n^2-1)], \qquad \text{and}  $$
$$ r = - \frac{4 (n-1)(n+3)}{\k^4} (\k^2 + 4c)^2 [\k^2 +  2c (n+1)] $$
we can verify that (97) holds by equating all components with zero. Note that the normal space to $\Bbb CQ^m$ in $H^{(1)}(m+1)$
at a point $P \in M^n $ is spanned by
vectors of the form $\tx , \s (X, Y), \s (X, \xi),$ $ \s (\xi , \xi), $ for $X, Y \in \Cal D.$
Moreover, the roots of the cubic equation $t^3 + p t^2 + q t + r = 0 $ are real and they are found to be
$$ \lambda_u = \frac{2c (n-1)(\k^2 + 4c)}{\k^2},
$$
$$ \lambda_v, \lambda_w = \frac{(\k^2 + 4c)[\k^2 + 4c (n+1)] \pm \sqrt{(\k^2 + 4c)[\k^6 - 12c \k^4 + 64c (n+1)^2]}}{2 \k^2}.
$$
Note that $\k^2 + 4c \neq 0 $ for a class$-B$ hypersurface. Equality of any two roots is possible only when $c = 1$ and $\lambda_u = \lambda_w $ (with minus sign at the radical), where $\k^2 = 2 (\sqrt{2m^2-1} - 1),$ identifying it as the 2-type example of Theorem 1 $(v).$ For the example given in Theorem 1 $(iv)$ we have $\k^2 = 4m $ and
$$ P(t) = [t^2 - \frac{4(m+1)(2m-1)}{m} t + 16 \frac{(m-1)(m+1)^2}{m} ][t - 8(m+1)],
$$
where the quadratic trinomial in the first pair of brackets is the minimal polynomial of that 2-type hypersurface according to (83), (84).
\qed
\enddemo
Now we are in a position to prove our classification result for 3-dimensional Hopf hypersurfaces of $\Bbb CQ^2(4c):$
\proclaim{Theorem 3} Let $M^3$ be a Hopf hypersurface of $\Bbb CP^2(4)$ with constant mean curvature.
Then $M^3$ is mass-symmetric and of 3-type in $H(3)$ if and only if $M^3$ is a class-$B$ hypersurface, that is, an open portion of a tube
of any radius $r \in (0, \pi /4)$ about the complex quadric $Q^1$
(equivalently the tube of radius $\frac{\pi}{4} - r$ about a canonically embedded $\Bbb RP^2$),
except when $\, \cot r = \sqrt 2 + \sqrt 3 $ and $\cot r = \sqrt{\sqrt 6 + \sqrt 7}.$
\endproclaim
\demo{Proof}
According to Corollary 2 such hypersurface has constant principal curvatures and thus it is one from the Takagi's list
in $\Bbb CP^2.$  Earlier analysis shows that it cannot be a  hypersurface of class $A_2,$ which is of 2-type when mass-symmetric
nor any of geodesic spheres, which are of 1- and 2-type. Standard examples of class $C, D,$ and $E$ need not be considered
because  of dimension restriction. Then Lemma 7 proves that class-$B$ examples are in fact the only ones, excluding the two
2-type tubes referred to in Lemma 4. \qed
\enddemo
\proclaim{Theorem 4} Let $M^3$ be a Hopf hypersurface of $\Bbb CH^2(-4)$ with constant mean curvature and
$(\text{tr}\, A)^2 \neq 4.$ Then $M^3$ is mass-symmetric and of 3-type in $H^1(3)$ if and only if $M^3$ is a
class-$B$ hypersurface, that is, an open portion of a tube of any radius $r > 0$ about a canonically embedded, totally real, totally geodesic $\Bbb RH^2 \subset \Bbb CH^2(4).$
\endproclaim
\demo{Proof}We know by Corollary 1 that the principal curvatures are constant and therefore examples are to be found among the standard 
ones from the Montiel's list. Every example of class $B$ is mass-symmetric and of 3-type by Lemma 7. Moreover, by (67), (69), 
with $c = -1, \; m = 2,$ we see that $\text{tr}\, A \neq 2$ for these hypersurfaces. A class$-A_0$ hypersurface in $\Bbb CH^m $ 
(a horosphere) is not of finite type since it satisfies $\Delta^2\tx = \text{const} \neq 0 $ [18]. Class$-A_1$ hypersurfaces (geodesic 
spheres and tubes about the complex hyperbolic hyperplane) are of 2-type. Class$-A_2$ hypersurfaces i.e. tubes about totally geodesic 
$\Bbb CH^k(-4), \; 1 \leq k \leq m-2, $ are of 3-type as shown before and among them there are some mass-symmetric ones. But because 
$m=2$ here, these examples are not possible. Note that an $A_2-$hypersurface degenerates into an $A_1-$hypersurface when $l = 0$ or $k = 0.$
\qed
\enddemo
The case of 3-type hypersurfaces of $\Bbb CH^2(-4)$ with $(\text{tr}\, A)^2 = 4$ is also interesting, but a different 
analysis is needed to study them. Because of this property, they are akin to the so-called Bryant surfaces in $\Bbb RH^3,$ see [17]. 
Also it would be interesting to determine the Chen-type of the standard examples of class $C, D$ and $E$ in $\Bbb CP^m(4)$ and, generally, 
study CMC Hopf hypersurfaces of 3-type in $\Bbb CQ^m(4c)$ when $m \geq 3.$ The techniques developed here can be modified to 
study curvature-adapted hypersurfaces of low type in quaternionic space forms and partly also in octonion planes, the topic that 
will be treated in our subsequent papers.

\enddocument